\documentclass[12pt,draftclsnofoot,onecolumn]{IEEEtran}

\usepackage[utf8]{inputenc} 
\usepackage[T1]{fontenc}
\usepackage{url}
\usepackage{ifthen}
\usepackage{cite}
\usepackage[cmex10]{amsmath} 
\usepackage{physics}
\usepackage{array}
\usepackage{makecell}
\usepackage{verbatim}
\usepackage{graphicx}
\usepackage{float}
\usepackage[colorinlistoftodos]{todonotes}
\usepackage{amssymb,amsfonts,bm}
\usepackage{amsthm}
\usepackage{tcolorbox}
\usepackage[framemethod=tikz]{mdframed}
\usepackage{mathtools,xparse}
\usepackage{tikz}
\usetikzlibrary{arrows}
\usetikzlibrary{shapes}
\usepackage{pgfplots}
\usetikzlibrary{plotmarks}
\pgfplotsset{compat=newest} 

\newtheorem{theorem}{Theorem}
\newtheorem{lemma}[theorem]{Lemma}
\newtheorem{corollary}[theorem]{Corollary}
\DeclareMathOperator*{\argmax}{arg\,max}

\def\cX{\mathcal{X}}

\def\bbeta{\bm{\beta}}

\def\pr{\text{Pr}}
\def\mgt{M^\textrm{GT}}

\def\sG{\text{sub-Gaussian}}

\def\sg{\text{sub-Gamma}}

\def\tv{\text{Var}}                       
\DeclarePairedDelimiter\ceil{\lceil}{\rceil}
\DeclarePairedDelimiter\floor{\lfloor}{\rfloor}

\def\fa{\{\alpha\}}
\def\ia{\lfloor\alpha\rfloor}
\def\Ia{\lceil\alpha\rceil}

\begin{document}
\title{Missing $g$-mass: Investigating the Missing Parts of Distributions}
\author{Prafulla~Chandra~and~Andrew~Thangaraj
  \thanks{P. Chandra and A. Thangaraj are with the Department of Electrical Engineering, Indian Institute of Technology Madras, Chennai, India 600036, Email: \{ee16d402, andrew\}@ee.iitm.ac.in. Parts of this paper appeared in the IEEE International Symposium on Information Theory 2019, Paris, France.}}

\maketitle

\begin{abstract}
Estimating the underlying distribution from \textit{iid} samples is a classical and important problem in statistics. When the alphabet size is large compared to number of samples, a portion of the distribution is highly likely to be unobserved or sparsely observed. The missing mass, defined as the sum of probabilities $\pr(x)$ over the missing letters $x$, and the Good-Turing estimator for missing mass have been important tools in large-alphabet distribution estimation. In this article, given a positive function $g$ from $[0,1]$ to the reals, the missing $g$-mass, defined as the sum of $g(\pr(x))$ over the missing letters $x$, is introduced and studied. The missing $g$-mass can be used to investigate the structure of the missing part of the distribution. Specific applications for special cases such as order-$\alpha$ missing mass ($g(p)=p^{\alpha}$) and the missing Shannon entropy ($g(p)=-p\log p$) include estimating distance from uniformity of the missing distribution and its partial estimation. Minimax estimation is studied for order-$\alpha$ missing mass for integer values of $\alpha$ and exact minimax convergence rates are obtained. Concentration is studied for a class of functions $g$ and specific results are derived for order-$\alpha$ missing mass and missing Shannon entropy. Sub-Gaussian tail bounds with near-optimal worst-case variance factors are derived. Two new notions of concentration, named strongly sub-Gamma and filtered sub-Gaussian concentration, are introduced and shown to result in right tail bounds that are better than those obtained from sub-Gaussian concentration. 
\end{abstract}

\textbf{Index terms-}  Missing mass, Good-Turing estimator, missing mass of a function, entropy, Mean squared error, minimax optimality, concentration, tail bounds, sub-Gaussian and sub-Gamma tails.   

\section{\textcolor{blue}{Introduction}}
Let $P$ be an arbitrary discrete distribution on an alphabet $\cX$. For $x\in\cX$, let $p_x\triangleq P(x)$. Let $X^n=(X_1,X_2,\ldots,X_n)$ be $n$ random samples with $X_i\sim P$ \emph{iid}.  Let $I(\cdot)$ and $E[\cdot]$ denote indicator random variables and expectations, respectively. For $x\in\cX$, let $N_x(X^n)=\sum_{i=1}^n I(X_i=x)$ denote the number of occurrences of $x$ in $X^n$. For a nonnegative integer $l$, let $\phi_l(X^n) \triangleq \sum_{x\in\cX} I(N_x(X^n)=l)$ denote the number of letters that occur $l$ times in $X^n$.

The problem of estimating $P$ from $X^n$, i.e. distribution estimation, is a classical problem in statistics. A common maxim in distribution estimation is that symbols that appear the same number of times in the sample are assigned the same probability in the estimate. To be more precise, for $l=0,1,\ldots$, let $S_l(X^n) \triangleq \{x\in\cX:N_x(X^n)=l\}$ denote the set of letters that have occurred $l$ times in $X^n$. In an estimate $\hat{P}$ of $P$, estimators studied in the literature typically enforce that $\hat{P}(u)=\hat{P}(v)$ if $u,v\in S_l$ (see, for example, \cite{Kamath15} and references therein). In the case of $l=0$, competitive distribution estimation against a "natural" genie that assigns equal probability to the missing letters $S_0(X^n)$ has been studied in \cite{Orlitsky2015}. 

Given the popularity of the above maxim, an interesting question to ask is whether we can test or infer if the probabilities of letters that have occurred the same number of times in the sample are likely to be all equal. In other words, is $P_l(X^n) \triangleq [p_x:x\in S_l(X^n)]$ close to being constant-valued or is it likely to be multi-valued? For $l=0$ and low $l$, this question is particularly interesting because it could reveal the structure of the missing part or sparsely observed part of the distribution. 

\subsection{Missing mass, Unseen species}
The missing mass, denoted $M_0(X^n,P)$ and defined as $M_0(X^n,P)=\sum_{x\in\cX}p_x I(N_x=0)$ is the sum of elements of $P_0(X^n)$ or total mass of unseen letters. The unseen species $\phi_0(X^n)$ is the number of unseen letters $|S_0(X^n)|$, and, in terms of estimation, is closely related to the support size $|\cX|$. Both these quantities, which involve the missing part of the distribution, have been studied extensively in prior work. 

Estimation of missing mass is widely used in domains like language modelling \cite{Gale1995,Chen1996} and ecology \cite{Chao1992}. A classical estimator for $M_0(X^n,P)$ is the Good-Turing estimator \cite{Good1953}, denoted $\mgt_0(X^n)$, defined as 
\begin{equation} 
\mgt_0(X^n) \triangleq \frac{\phi_1(X^n)}{n},
\end{equation}
where $\phi_1(X^n)$ is the number of letters that have occurred once in $X^n$. The problem of estimation and concentration of missing mass and the properties of the classical Good-Turing estimator have been studied in \cite{McAllester2000, Orlitsky2003, mansour2009multiple, Wagner2006, wagner2007better, ohannessian2012rare, Acharya13, Orlitsky2015}. Missing mass has been used to lower bound the total variation distance of the empirical distribution of i.i.d samples to their original distribution in \cite{kontorovich_2020}.

While missing mass estimation under squared-error loss is consistent over all distributions with no further restrictions, consistent estimation of support size or unseen species requires a restriction to the class of distributions in which $p_x\notin(0,1/k)$, where $k$ is a parameter. This problem has been studied by several authors, and optimal minimax results were provided in  \cite{WuYangChebyshev2019} (see references here for other works). 

\subsection{Proximity of $P_0(X^n)$ to uniformity and missing $g$-mass}
An interesting question at this juncture is whether anything more can be inferred about $P_0(X^n)$ beyond its sum and size. Let us consider estimating the closeness to uniformity of $P_0(X^n)$, which could be measured as a distance between $P_0(X^n)=[p_x:N_x(X^n)=0]$ and $[\delta\ \delta\ \cdots\ \delta]$ (of the same length), for instance, in the following two possible ways:
\begin{enumerate}
    \item Relative log distance:
    \begin{align}
        \sum_{x \in \cX: N_x=0} p_x \log_2 (p_x/\delta)  =   M_0\log_2(1/\delta) - \sum_{x\in\cX} (-p_x\log_2 p_x)\,I(N_x=0). \label{eq:KL_div_missing}
    \end{align}
    \item Squared-error distance: 
    \begin{align}
        \sum_{x \in \cX: N_x = 0} (p_x - \delta)^2= \delta(\delta\phi_0-2M_0) + \sum_{x \in \cX} p_x^2\,I(N_x=0). \label{eq:L2_loss_missing}
    \end{align}
\end{enumerate}
where $M_0, \phi_0$ and $N_x $ are compact notation for $M_0(X^n,P), \phi_0(X^n)$ and $N_x(X^n)$, respectively. In both of the above distance expressions, there appears a summation term of the form 
\begin{equation}
M_{0,g} (X^n, P) \triangleq \sum_{x \in \mathcal{X} } g(p_x)\, I(N_x(X^n) = 0),
\label{eq:unseen_def}
\end{equation} 
where the function $g(p)$ is $-p\log_2 p$ in \eqref{eq:KL_div_missing} and $p^2$ in \eqref{eq:L2_loss_missing}. We refer to the quantity $M_{0,g}(X^n,P)$ defined in \eqref{eq:unseen_def} as the missing $g$-mass. We see that missing $g$-mass for some specific $g$ can possibly reveal more information about the structure of the missing part of the distribution $P_0(X^n)$. We further illustrate this with simulations in Section \ref{subsec:unseen_part_sim}. 

The two special cases in \eqref{eq:KL_div_missing} and \eqref{eq:L2_loss_missing} are studied in more detail in this article, and are defined separately below. 
\begin{equation}
   M_{0,\alpha}(X^n,P) \triangleq \sum_{x\in\cX} p_x^{\alpha} \, I(N_x(X^n)=0) ,\ \alpha > 0, \label{eq:M0_alpha_def} \end{equation}
is called the missing mass of order-$\alpha$, 
and 
\begin{equation}
   H_0(X^n,P) \triangleq \sum_{x\in\cX} p_x \ \log_2 (1/p_x) \, I(N_x(X^n)=0) \label{eq:H0_def} 
   \end{equation}
 is the missing Shannon entropy of $X^n$. We will drop the arguments $X^n$ and $P$ whenever possible.


The problem of estimating additive functions $\sum_{x\in\cX}g(p_x)$ of the distribution $P$ from $n$ random samples has been studied before \cite{jiao2015minimax}. Estimation of Shannon entropy, which is an additive function with $g(p)=-p\log_2 p$, and R\'{e}nyi entropy of order-$\alpha$, which contains the additive function with $g(p)=p^{\alpha}$, have been of particular interest \cite{Acharya17Renyi}. While the applications of estimating Shannon entropy are numerous, estimation of R\'{e}nyi entropy has applications, for instance, in the problem of guessing \cite{arikan1996inequality, hanawal2011guessing}. 

To the best of our knowledge, this is the first work to consider the above generalization of missing mass to missing $g$-mass. 

\subsection{Minimax estimation setting}

Our study of estimation will be in the worst case minimax sense. An estimator $\widehat{G}_0(X^n)$ for $M_{0,g}(X^n,P)$ is a mapping from $\cX^n$ to $\mathbb{R}$. For a distribution $P$, the $L_2^2$ or squared-error risk of the estimator $\widehat{G}_0(X^n)$  is 
\begin{equation}
    R_{n,g}(\widehat{G}_0, P) \triangleq E_{X^n \sim P} [ (\widehat{G}_0(X^n) - M_{0,g} (X^n, P) )^2]. \label{eq:risk}
\end{equation}
The worst case risk of $\widehat{G}_0(X^n)$ is 
\begin{equation}
    R_{n,g}(\widehat{G}_0) \triangleq \max_P R_{n,g}(\widehat{G}_0, P). \label{eq:max_risk}
\end{equation}
The minimax risk of estimating $M_{0,g} (X^n, P)$ is 
\begin{equation}
    R_{n,g}^* \triangleq \min_{\widehat{G}_0 } R_{n,g}(\widehat{G}_0). \label{eq:minimax_risk}
\end{equation}
As is standard, the goal of our study is to characterize the rate of decay of $R_{n,g}^*$ with $n$ for a given function $g(p)$. For the estimation part, in this article, we will primarily consider $g(p)=p^{\alpha}$, where $\alpha=1,2,\ldots$. In Section \ref{subsec:Minimax_est}, we introduce the generalised Good-Turing estimators that consistently estimate $M_{0,\alpha},$ for $\alpha \in \mathbb{N} \triangleq \{1,2,\ldots\}$ and show by simulations in Section \ref{subsec:unseen_part_sim} that using these estimators lead to reliable inference about the uniformity of the missing probabilities. We also present simulations that use the estimates of $M_{0,\alpha}$ for investigating the values of missing probabilities along with their multiplicities.


We make an important remark regarding the choice of squared-error loss function. In \cite{ohannessian2012rare,Ohannessian19}, it has been shown that the missing mass $M_0$ cannot be learned uniformly over all distributions when the loss function is the relative error (i.e. $|\frac{\widehat{M_{0}}}{M_{0}} -1|$). In regimes where $M_0\to0$ as $n\to\infty$, the use of relative error is, therefore, important. 
However, in the large alphabet regime where $|\cX|$ is expected to be larger than $n$, $M_0$ will take on a high value and be a significant fraction with nonzero probability. As a simple example, when $cn$ letters have probability in the interval $[c'/n,c''/n]$, $E[M_{0}]=\sum_{x}p_x(1-p_x)^n > (cn)(c'/n)(1-c''/n)^n > cc'(1-c'')$. Given the concentration properties of $M_{0}$ \cite{berend2013}, we get that $M_{0}$ is bounded away from 0 with high probability as well in such distributions. Several variations of such examples occur in the large alphabet regime, where $|\cX|$ is larger than $n$. Since estimation of missing mass is typically important only when it takes on a significant value, the squared-error loss has been popular in the study of missing mass, and we adopt the same loss function for missing $g$-mass as well.



\subsection{\textcolor{blue}{Novel Concentration Results}} 
The estimation of missing mass is considered to be complicated primarily because missing mass is a random variable dependant on the samples $X^n$ and the distribution $P,$ while typical quantities for estimation are functions of $P$ alone. However, missing mass is known to concentrate about its expected value \cite{McAllester2000,berend2013,ben-hamou2017}, and the concentration phenomenon plays an important role in the success of estimation of missing mass.  

Most concentration results involve bounding the moment generating function (MGF) of a random variable of interest, say $Z$, with a known closed form MGF of a random variable $Y$ (say, Gaussian) having a suitable variance $v$. A highly desirable concentration result is sub-Gaussianity with $v = \text{Var}(Z)$. However, there exist asymmetric random variables like the missing mass $M_0$, which is closer to $0$ than $1$, 
and has a left tail decaying faster than the right tail. For such a random variable $Z$, the best sub-Gaussian bound for the lighter tail might have a variance factor $v > \text{Var}(Z)$ leaving the problem of bounding the MGF of $Z$ with that of a suitable (non-Gaussian) $Y$ with $v = \text{Var}(Z)$ as an open question. To answer this question and thereby achieve improved right tail bounds for $M_0,$ we propose two new types of concentration called the \textit{Strongly sub-Gamma} and the \textit{filtered sub-Gaussian} concentrations. 

\noindent\textbf{Strongly sub-Gamma:} We say a random variable $Z$ is strongly sub-Gamma with variance factor $v$ and scale parameter $c$ if 
\begin{equation}
    \ln E[\exp{\lambda (Z - E[Z])}] \ \leq \ \frac{v}{c^2} \ \ln \ \frac{e^{-\lambda c}}{1 - c\lambda}, \ \lambda \in (0,1/c). \label{eq:strongly-sub-gamma-MGF_ub}
\end{equation}
Using \eqref{eq:strongly-sub-gamma-MGF_ub} in the Chernoff method and simplifying, we get the following right tail bound:
\begin{equation}
    \pr(Z\ge E[Z]+\epsilon)\le \exp\left\{ -\frac{\epsilon^2}{2v} \left( 1 - ({2}/{3}) \frac{c\epsilon}{v} \right) \right\}  \label{eq:strongly-sub-gamma-tail_ub}
\end{equation}

\noindent\textbf{Filtered sub-Gaussian:} We say a random variable $Z$ has filtered sub-Gaussian concentration with variance factor $a_2$ and filter $h(\lambda)$ if 
\begin{equation}
    \ln E[\exp{\lambda (Z - E[Z])}] \ \leq \ \frac{\lambda^2a_2}{2} +  h(\lambda) \label{eq:filt-sub-gaussian-MGF_ub}
\end{equation}
If $Z$ satisfies \eqref{eq:filt-sub-gaussian-MGF_ub} with $h(\lambda) = \frac{v}{c^2} \ln \ \frac{e^{-\lambda c}}{1 - c\lambda}$ for $\lambda \in (0,1/c)$ with $v \geq a_2,$ we get
\begin{equation}
    \pr(Z\ge E[Z]+\epsilon)\le \exp\left\{ -\frac{\epsilon^2}{2(v+a_2)} \left( 1 - ({2}/{3}) \frac{cv\epsilon}{(a_2 + v)^2} \right) \right\} \label{eq:poly-strongly-sub-gamma-tail_ub}
\end{equation}
using the Chernoff method. Note that the tail bounds in \eqref{eq:strongly-sub-gamma-tail_ub} and \eqref{eq:poly-strongly-sub-gamma-tail_ub} have only one additional term inside the exponent when compared to the right tail bounds obtained from sub-Gaussianity with variance factors $v$ and $a_2+v$, respectively. If the best possible sub-Gaussian tail bound for $Z$ has variance $> \text{Var}(Z),$ then characterising $Z$ to be Strongly sub-Gamma with variance factor $v = \text{Var}(Z)$ (or Filtered sub-Gaussian with $a_2 + v = \text{Var}(Z))$, results in a right tail bound smaller than the best possible sub-Gaussian bound (and much closer to the desired sub-Gaussian bound with $v= \text{Var}(Z))$ over a restricted range of $\epsilon$.
In our results described in Section \ref{sec:Our_results}, we show that $M_0$ is uniformly \textit{Strongly sub-Gamma} and \textit{filtered sub-Gaussian} with variance factor $v \approx \max_{P} \ \text{Var}(M_0)$, while the best possible uniform sub-Guassian bound has variance $= 0.5/n >  \max_{P} \ \text{Var}(M_0)$. 
Thus, we improve the best possible uniform sub-Gaussian right tail bound for $M_0$ over a restricted range  of $\epsilon$.

We also extend these characterizations to the concentration of $M_{0,g}(X^n,P)$ 
for two classes of functions $g(p)$, and provide specialised results for the case $g(p)=p^{\alpha}$ and $g(p)=-p \log_2 p$. Our approach for proving concentration is different from previous ones, and, arguably, provides one of the simplest proofs for concentration results of missing mass. 

The rest of this article is organised as follows. 
We present a summary of our main results together with a brief description of the relevant prior results from the literature in Section \ref{sec:Our_results}\textcolor{blue}{, and present our experiments on estimating the missing part of the distribution $P$  and testing its closeness to uniformity}. In Section \ref{sec:minimax_M_0_alpha_proof}, we provide the proofs for our results on the minimax estimation of $M_{0,\alpha}$ for $\alpha\ge 1$ under $L^2_2$ risk. In Section \ref{sec:conc_proofs}, we present the proofs for our concentration and tail bound results for $M_{0,g}$. Section \ref{sec:rec_proof} presents a proof for a lemma outlining the choice of parameters in the concentration results. Section \ref{sec:corollary_proof} has derivations of specific, simplified tail bounds 
for generalized missing mass and missing Shannon entropy. We conclude in Section \ref{sec:conclusion} with some remarks. Technical aspects of some proofs are collected in the appendix.

\section{Summary of results} \label{sec:Our_results}
In this section, we provide details of the problems addressed in this work, summarize the main results and discuss them. 

\subsection{Minimax estimation results for $M_{0,\alpha}$} \label{subsec:Minimax_est}
Let $R^*_{n,\alpha}\triangleq R^*_{n,p^\alpha}$ denote the minimax risk of estimating the order-$\alpha$ missing mass $M_{0,\alpha}$ under $L_2^2$ loss as defined in \eqref{eq:minimax_risk}. Let $\mathbb{N}=\{1,2,\ldots\}$ denote the set of positive integers. For two non-negative sequences $a_n$, $b_n$, the notation $a_n =_n b_n$ denotes that $a_n = b_n \pm o(b_n).$ The notations $a_n \le_n b_n$ and $a_n \ge_n b_n$ are similarly defined.  

As defined above, $R^*_{n,p}$ denotes the minimax risk of estimating the missing mass $M_0=\sum_x p_x\,I(N_x=0)$ under $L_2^2$ loss. While $R^*_{n,p}$ was known to be $O(1/n)$ on the basis of concentration results \cite{McAllester2000}, bounds with small constants were found in \cite{Rajaraman17}, which showed that 
\begin{align} 
\frac{0.25}{n} & \leq_n R_{n,p}^* \leq_n \frac{0.6179}{n}. \label{eq:minmax_M0_1} 
\end{align} 
In \cite{Acharya17}, the bounds were further improved to 
\begin{align}
    \frac{0.570}{n} &\leq_n R_{n,p}^* \leq_n \frac{0.608}{n}. \label{eq:minmax_M0_2}
\end{align} 
The upper bound in the above results comes from the worst case risk of the Good-Turing estimator. In our first result, presented in the following theorem, we provide matching upper and lower bounds (in their dependance on $n)$ on $R^*_{n,\alpha},$ for $\alpha \in \mathbb{N.}$   
\begin{theorem} 
 For a positive integer $\alpha\in \mathbb{N}$ and $n > 2\alpha$, 
\begin{align}
     \frac{c_l}{n^{2\alpha -1}} \le_n R_{n,\alpha}^* \le_n \frac{c_u}{n^{2\alpha -1}},
    \label{eq:minimax_M_0_alpha_int}
\end{align}
where $c_l$, $c_u$ are positive constants.
\label{thm_part:minimax_M_0_alpha_int}


\label{thm:minimax_M_0_alpha}

\end{theorem}
Paraphrasing the theorem, the minimax risk of the order-$\alpha$ missing mass for a positive integer $\alpha$ falls as $1/n^{2\alpha-1}$. In a proof presented in Section \ref{sec:minimax_M_0_alpha_proof}, we show the lower bound in \eqref{eq:minimax_M_0_alpha_int} using Dirichlet priors and the upper bound using a generalized Good-Turing estimator for $M_{0,\alpha}(X^n,P)$, denoted $\mgt_{0,\alpha}(X^n)$, defined as 
\begin{equation}
    \mgt_{0,\alpha}(X^n)=\dfrac{\phi_{\alpha}}{\binom{n}{\alpha}},\ \alpha\in\mathbb{N}.
    \label{eq:mgt}
\end{equation}
To understand why the above estimator works well, let us consider its bias, which can be simplified as follows:
\begin{align}
E[\mgt_{0,\alpha}(X^n)-M_{0,\alpha}(X^n,P)]&=\sum_{x\in\cX}\dfrac{1}{\binom{n}{\alpha}}\,\pr(N_x=\alpha)\, -\, p^\alpha_x\,\pr(N_x=0)\nonumber\\
&= \sum_{x\in\cX} p^\alpha_x(1-p_x)^{n-\alpha}\left[1-(1-p_x)^\alpha\right] \quad (\ge 0)\nonumber\\
|E[\mgt_{0,\alpha}(X^n)-M_{0,\alpha}(X^n,P)]|&\overset{(a)}{\le} \dfrac{\alpha}{\binom{n}{\alpha}} \sum_{x\in\cX} p_x\left[ \binom{n}{\alpha}p^\alpha_x(1-p_x)^{n-\alpha} \right]\label{eq:pxterm}\\
&\overset{(b)}{\le} \dfrac{\alpha^{\alpha+1}}{n^\alpha},
\end{align} 
where the change to absolute value in the LHS of $(a)$ follows because the previous expression is clearly non-negative, and the RHS of $(a)$ follows from the inequality $1-(1-p_x)^\alpha\le \alpha p_x$ ($\alpha\ge 1$). Finally, $(b)$ follows because $\binom{n}{\alpha}p^\alpha_x(1-p_x)^{n-\alpha}\le 1$ and $\binom{n}{\alpha}\ge (n/\alpha)^\alpha$. Hence, the bias falls as $1/n^{\alpha}$. The squared-error risk computation is more involved, and, as shown in the proof in Section \ref{sec:minimax_M_0_alpha_proof}, it falls as $1/n^{2\alpha-1}$ matching the lower bound.

At this juncture, we draw attention to \eqref{eq:pxterm}, where the $p_x$ term inside the summation results in an averaging of a bounded function from $[0,1]\to\mathbb{R}$. The average is upper bounded by the maximum of the function, and this bound holds for all distributions with no assumptions needed on the alphabet size. This trick is a recurring one in the analysis of estimation of missing mass and its variants. Notice that the division by $\binom{n}{\alpha}$ in the definition of the estimator $\mgt_{0,\alpha}$ in \eqref{eq:mgt} is instrumental in the additional $p_x$ term appearing inside the summation in the bias expression. \textcolor{blue}{In this work, we follow the traditional approach to studying missing mass using the above trick and do not consider the use of the support size as a parameter.} 

\textcolor{blue}{The estimand $M_{0,\alpha}$ has an expected value $E[M_{0,\alpha}]=\sum_x p_x^{\alpha}(1-p_x)^n\ge (c)^{\alpha}(1-c')/n^{\alpha-\epsilon}$ for distributions with at least $n^\delta$ letters with probabilities in the range $\left[\dfrac{c}{n^{1+(\delta-\epsilon)/\alpha}},\dfrac{c'}{n}\right]$, $\delta>\epsilon$. The concentration results for $M_{0,\alpha}$ (see Section \ref{sec:conc_results}) imply that $M_{0,\alpha}\ge c''/n^{\alpha-\epsilon+o(1)}$ with high probability. When compared to the bias of $1/n^{\alpha}$, we see that the estimand is larger for the considered distributions, which are quite common in the large alphabet regime. However, amongst all distributions, there are likely to be some light-tailed distributions, where $M_{0,\alpha}$ decays as $1/n^{\alpha}$ or faster. In such light-tailed cases, relative error should be considered and it is likely an impossibility result similar to that of \cite{ohannessian2012rare}\cite{Ohannessian19} holds true. However, the estimation of the missing structure is not as interesting in such light-tailed regimes because the missing part of the distribution does not become substantial enough to cause any effect on distribution estimation. In this work, we consider squared error loss, and relative error could be a topic for future study.}   

\noindent\emph{Remarks on generalization}: (1) The bias upper bound of $O(1/n^{\alpha})$ above can be generalized readily to the missing mass of $g(p)$ that satisfies
\begin{equation}
    \frac{g(p)(1-p)^{\alpha}}{p^{\alpha}} \in [a - O(p\,h(p)), a + O(p\,h(p))] 
\end{equation}
$\alpha\in\mathbb{N}$, $a$ is a constant and $\binom{n}{\alpha}p^{\alpha}(1-p)^{n-\alpha}h(p)\le O(1)$ by using the Good-Turing estimator $a\mgt_{0,\alpha}$. Linear combinations of a constant (possibly, $o(n)$) number of such $g(p)$'s can also be estimated with diminishing bias by using corresponding linear combinations of the Good-Turing estimators. Examples include $g(p)=p^c(1-p)^d$ ($c\in\mathbb{N}$, $d>1$), $p^ce^{-dp}$ ($c\in\mathbb{N}$) and so on. 

\textcolor{blue}{(2) Polynomial approximation \cite{WuYangChebyshev2019} can be readily used to extend estimators of $M_{0,\alpha}$, $\alpha\in\mathbb{N}$, to obtain estimators of $M_{0,g}$ for functions $g(p)$ that are sufficiently smooth and bounded for $p\in[0,1]$. Since this is a straight-forward use of the method described in \cite{WuYangChebyshev2019}, we have not elaborated on it in further detail. However, for interesting functions such as $g_1(p)=-p\log p$ and $g_2(p)=p^{\alpha}$, $\alpha\notin\mathbb{N}$, because $g'_1(p), g_2^{(\lceil\alpha\rceil)}(p)\to\infty$ as $p\to0$, a support size assumption on the distribution is required to ensure that approximation is needed only on an interval that is away from the origin. In this initial work, we have not considered support size as a parameter, and mainly focus on the case $g(p)=p^{\alpha}$, where $\alpha$ is a positive integer. Study of general classes of $g(p)$ and modifications required in polynomial approximation methods may be taken up in future work.}

\textcolor{blue}{(3) Finally, generalizing from missing mass $M_{0,g}$ to  $M_{l,g},$ the total mass of letters that have appeared $l$ times is interesting to consider. We allude to this generalization in the experimental results, and remark that most of our results for $M_{0,g}$ will generalise to $M_{l,g}$ with simple changes.}

\subsection{ \textcolor{blue}{Experiments}}\label{subsec:unseen_part_sim}
We illustrate the use of missing $g$-mass by describing two techniques to infer about the structure of the vector $P_0 = \{p_x : N_x = 0 \}$ using estimation of $M_{0,\alpha}$ with $\alpha \in \mathbb{N}.$
\subsubsection{Tests for uniformity of $P_l$} 
One approach to test for uniformity of $[p_x: x \in S_l(X^n)]$ is to consider integer power sums of $p_x,$ for $x \in S_l,$ defined as follows:
\begin{equation}
    M_{l,k}(X^n,P)\triangleq\sum_{x\in\cX} p_x^k\, I(N_x(X^n)=l) = \sum_{x\in S_l(X^n)} p_x^k.
\end{equation}
For $l=0,$ the above summation is the missing $g$-mass $M_{0,g}(X^n,P)$ with $g(p)=p^k$.

To reduce clutter, we will drop the arguments $X^n$ and $P$ appearing in the defined quantities. If for all $x\in S_l$, $p_x\approx p_l$, we expect that $M_{l,k}\approx |S_l|p_l^k$. So, if $[p_x:x\in S_l(X^n)]$ is close to constant-valued, we expect that the ratios $M_{l,k-1}/M_{l,k}\approx1/p_l$ and $(M_{l,1}/M_{l,k})^{1/(k-1)}\approx1/p_l$ for $k=2,3,\ldots$. This is further justified by the following versions of Cauchy-Schwarz and H\"{o}lder's inequalities for positive integers $a_1,\ldots,a_m$: 
\begin{align}
    \frac{a_1^k+\cdots+a_m^k}{a_1^{k+1}+\cdots+a_m^{k+1}}\le \frac{a_1^{k-1}+\cdots+a_m^{k-1}}{a_1^k+\cdots+a_m^k}\le\cdots&\le\frac{m}{a_1+\cdots+a_m},\\
    \left(\frac{a_1+\cdots+a_m}{a_1^k+\cdots+a_m^k}\right)^{\frac{1}{k-1}}&\le \frac{m}{a_1+\cdots+a_m}.
\end{align}
Equality holds in the above if and only if $a_1=\cdots=a_m$. 

\begin{figure}[htb]
    \centering
    \includegraphics[width=7in]{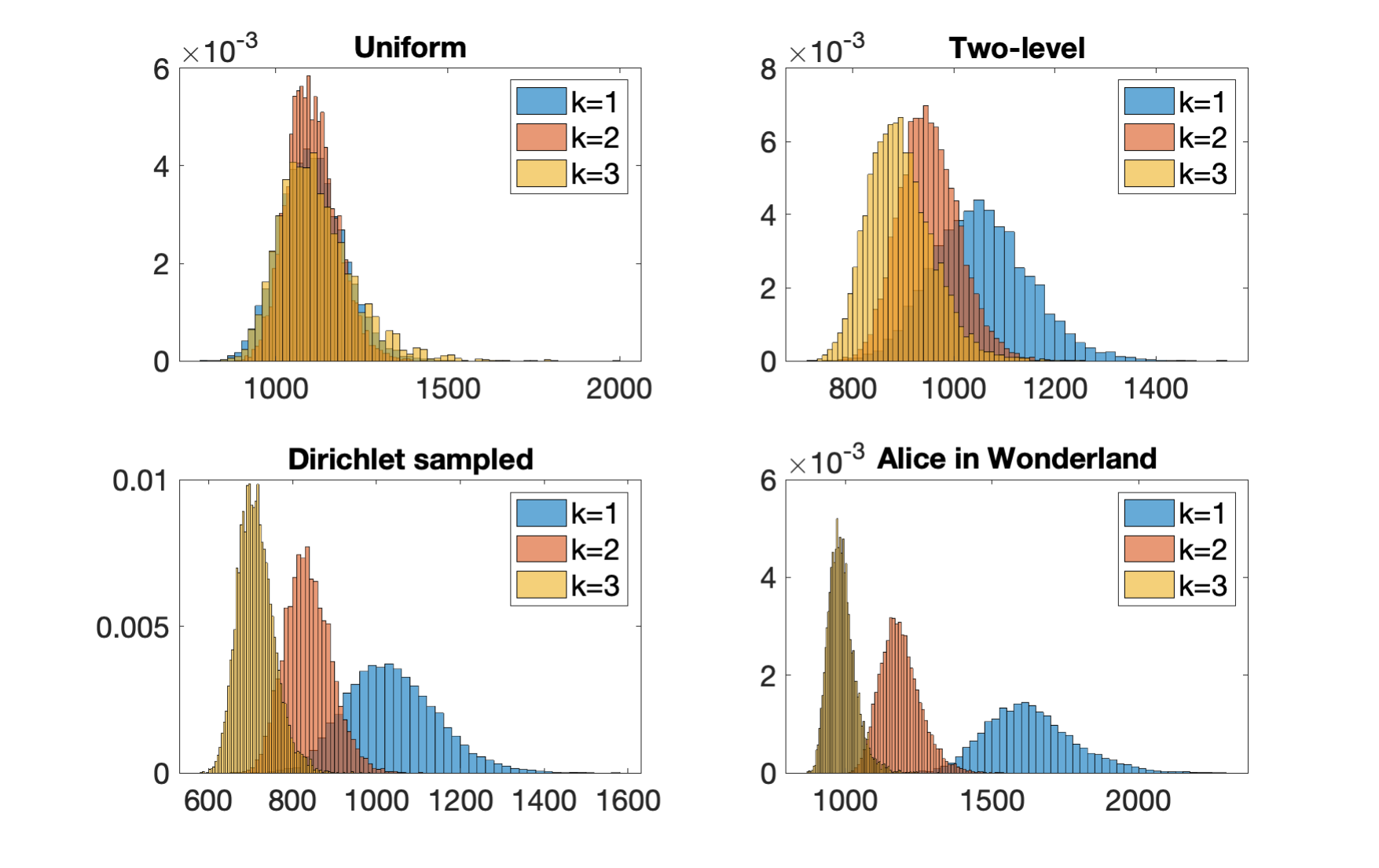}
    \caption{Histograms of candidate test statistics $(\mgt_{0,1}/\mgt_{0,k})^{1/(k-1)}$ for uniformity for $k=2,3,4$. $x$-axis is the value of the test statistic, and $y$-axis is the frequency.}
    \label{fig:hist}
\end{figure}

To illustrate the above, we use the generalized Good-Turing type estimator $\mgt_{0,k}(X^n)$ for $M_{0,k}.$ Fixing $n=1000$, we compute and plot histograms of $(\mgt_{0,1}/\mgt_{0,k})^{1/(k-1)}$ over $10000$ trials for 4 different distributions - (i) Uniform over 1100 letters, (ii) Two-level: 550 letters with probability 3/2200 and 550 letters with probability 1/2200, (iii) Distribution on 1100 letters sampled randomly from Dirichlet$(1/1100,\ldots,1/1100)$, (iv) Distribution of 2889 words from the text of the novel "Alice in Wonderland". The plot is shown in Fig. \ref{fig:hist}.
From the figure, one can see that the histograms are overlapping significantly for the uniform distribution and are overlapping partially for the two-valued distribution. The histograms separate out for the Dirichlet-sampled and ``Alice in Wonderland'' distributions, which show a significant level of nonuniformity. A reasonable conclusion from the plot is that the values of $(\mgt_{0,1}/\mgt_{0,k})^{1/(k-1)}$ for $k=2,3,4$ appear to act as useful test statistics for uniformity of the probability of missing letters.

\subsubsection{Estimating rare probabilities and multiplicities}
\begin{figure}[htb]
    \centering
    \includegraphics[width=7in]{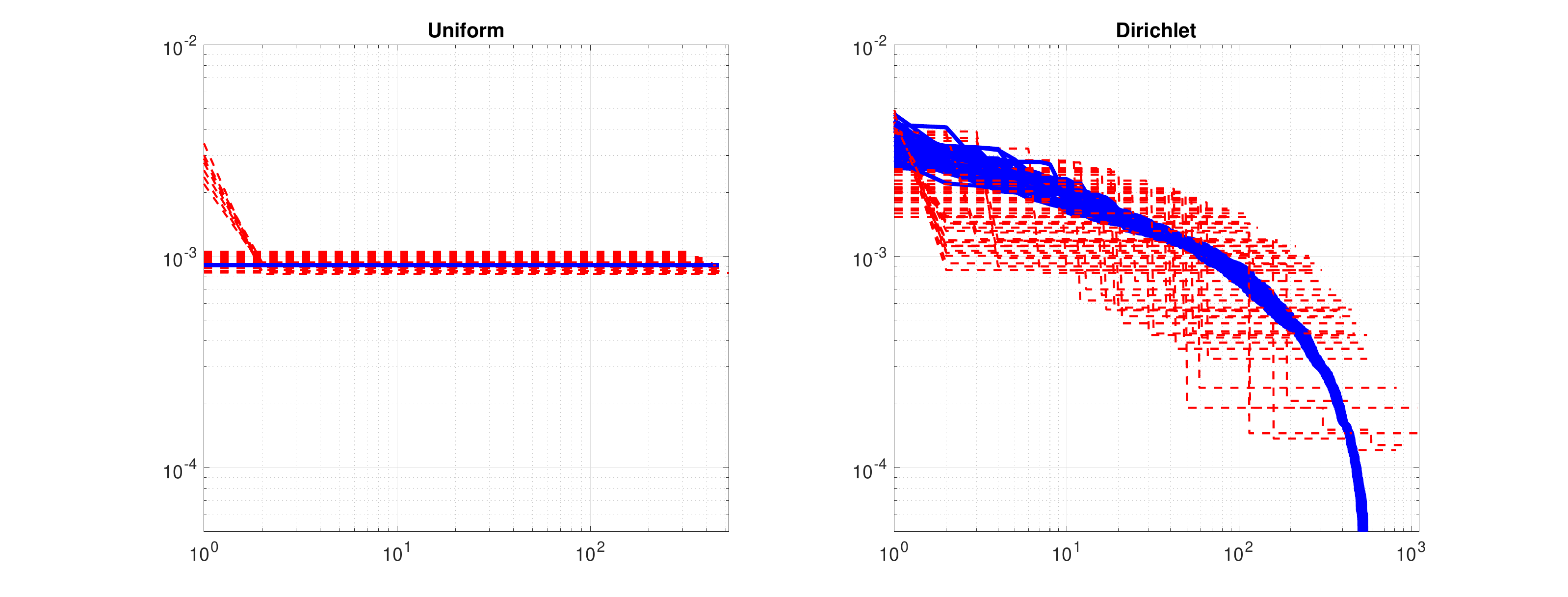}
    \caption{Actual and estimated distribution of missing letters (in descending order). Blue solid: actual, Red dashed: estimate.}
    \label{fig:est0}
\end{figure}

Another interesting objective is to estimate the vector $[p_x:x\in S_l(X^n)]$ (up to ordering) using estimates of the integer power sums \begin{equation}
    M_{l,k}(X^n,P)\triangleq\sum_{x\in\cX} p_x^k\, I(N_x(X^n)=l) = \sum_{x\in S_l(X^n)} p_x^k.
\end{equation}
A Good-Turing type estimate for $M_{l,k}$ can be derived to be 
\begin{equation}
    \mgt_{l,k}(X^n)=\frac{\binom{n}{l}\phi_{k+l}(X^n)}{\binom{n}{k+l}}.
\end{equation}
Estimators of this form for $k =1$ have already been studied and used for distribution estimation in the literature \cite{McAllester2000,Orlitsky2015, hao_doubly_comp19,Painsky21}. \textcolor{blue}{By following a method similar to the proof of \eqref{eq:minimax_M_0_alpha_int}, the MSE of $\mgt_{l,k}(X^n)$ can be shown to be of the order $O(1/n^{2k-1}).$}
As before, we fix $n=1000$ and consider the uniform distribution on 1100 letters and a distribution sampled from Dirichlet$(1/1100,\ldots,1/1100)$ for comparison. For $l=0,1$, samples from the above two distributions yield reasonable estimates of $M_{l,k}$ for $k=1,2,3,4,5$. Working with estimates of 5 power sums, we will estimate $[p_x:x\in S_l(X^n)]$ (in, say, descending order) as a vector of 5 different values $p_1,\ldots,p_5$ each occurring with respective multiplicities $m_1,\ldots,m_5$. This estimation could be accomplished, for instance, by solving the following optimisation problem:
\begin{equation}
    \min_{p_i,m_i} \sum_{j=1}^5 (m_1p_1^j+\cdots+m_5p_5^j - \mgt_{l,j})^2.
\end{equation}
By relaxing $m_i$ to be real, the above problem can be readily solved using several computational tools. For $l=0$, Fig. \ref{fig:est0} shows plots of $[p_x:x\in S_l(X^n)]$ and its estimate in descending order for the two chosen distributions for 50 trials.
We observe that the estimate is close to the actual distributions in all of the cases.

\begin{figure}[htb]
    \centering
    \includegraphics[width=7in]{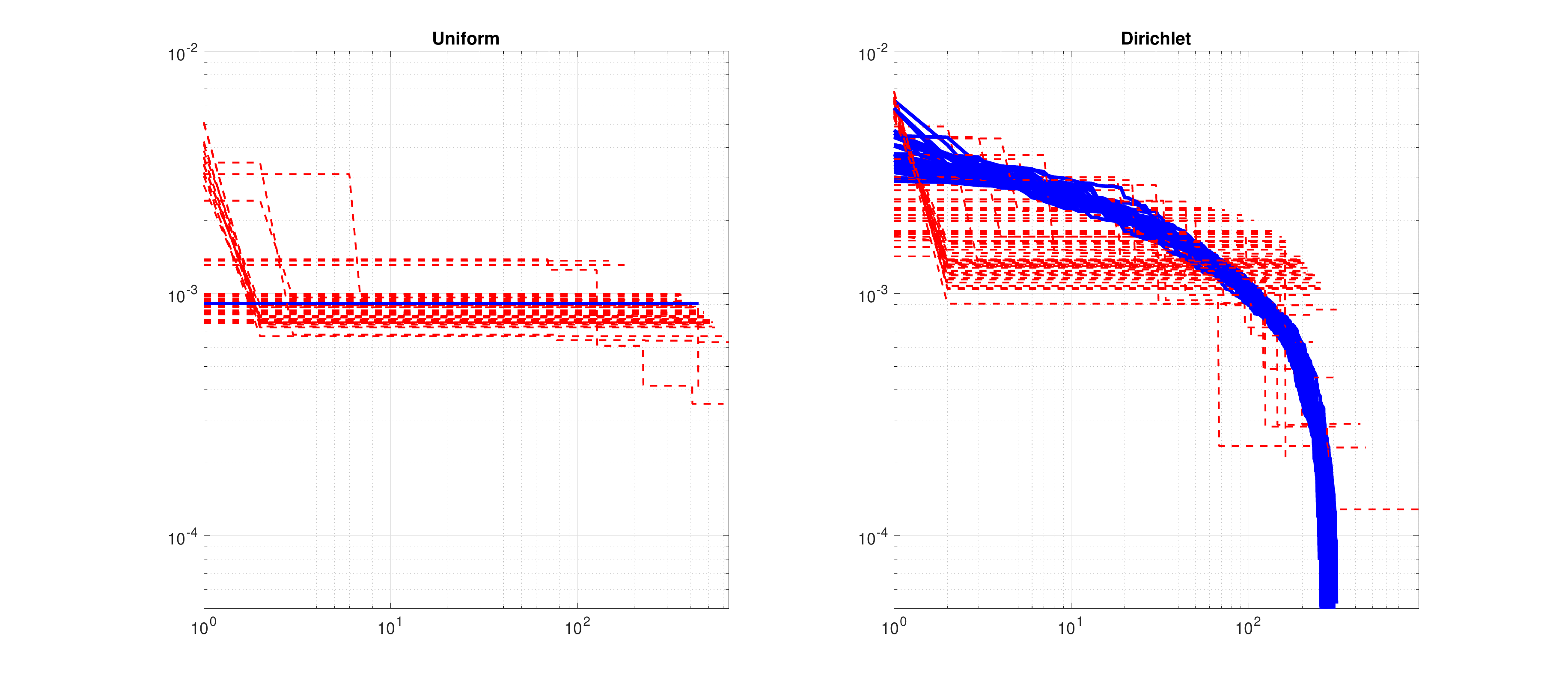}
    \caption{Actual and estimated distribution of letters that appear once (in descending order). Blue solid: actual, Red dashed: estimate.}
    \label{fig:est1}
\end{figure}

A similar plot for $l=1$ is shown in Fig. \ref{fig:est1}.
As observed, the estimates vary more than the $l=0$ case, and generally appear to be weaker than the results for $l=0$. Therefore, it is important to theoretically study the limits of estimation of $[p_x:x\in S_l(X^n)]$ and associated quantities.

In summary, the generalization of missing mass that we have proposed offers the possibility of exploring the structure of the missing distribution and the sparsely observed distribution. We have demonstrated two explicit applications that test for uniformity of missing letters and attempt to estimate the distribution of missing letters and letters that occur once, up to ordering.
\textcolor{blue}{ These testing and estimation problems introduced above aim to infer more about the probabilities of the letters that were not seen in the samples, departing from the already well-studied measures like missing mass and the number of unseen letters (or species) that do a coarse characterization of the unseen part of $P.$ Inferring the finer details about the missing probabilities, like the proximity of the missing part of the $P$ vector to uniformity or the individual probabilities of the missing letters, could lead to better estimation of $P$ from $X^n.$} 

\subsection{Concentration results and new tail bounds}
\subsubsection{Prior Work} 
\textcolor{blue}{Various concentration results (uniform over all distributions $P$) and tail bounds have been shown for missing mass $M_0(X^n,P)$ } \cite{McAllester2000, McAllester2003, berend2013, ben-hamou2017, Chandra_ncc}.
We quickly recall the general setting for such concentration results. For a random variable $Z$, let 
\begin{equation}
L_Z(\lambda)\triangleq \log(E[e^{\lambda(Z-E[Z])}])    
\end{equation}
denote the log Moment Generating Function (MGF). A random variable $Z$ is said to be $\sG$ with variance factor $v$, denoted $Z\sim\sG(v)$, on the right tail if $L_Z(\lambda)\le \lambda^2v/2$, $\lambda > 0$. If $Z$ is $\sG(v)$ on the right tail, by the standard Chernoff method, we have the right tail bound
\begin{equation}
    \pr(Z\ge E[Z]+\epsilon)\le e^{-\epsilon^2/2v},\ \epsilon \ge 0.
\end{equation}
For sub-Gaussian left tail bounds, we will need ${-}Z\sim\sG(v)$ on the right tail, or equivalently $L_Z(\lambda)\le \lambda^2v/2$, $\lambda < 0$.

One of the first concentration results for missing mass was shown in \cite{McAllester2000}, where a high probability bound on a letter not being seen in $X^n$ is used in McDiarmid's inequality \cite{boucheron2013} to show the right tail bound $\pr(M_0 \geq E[M_0]+\epsilon)\le e^{-n\epsilon^2/3}$, $\epsilon \ge 0$. In \cite{McAllester2003}, the above bound was improved to 
\begin{equation}
    \pr(M_0\ge E[M_0]+\epsilon) \le e^{-n \epsilon^2},\ \epsilon \ge 0,
    \label{eq:right_tail_best_sG}
\end{equation}
by showing that $M_0(X^n,P)\sim\sG(0.5/n),$ \textcolor{blue}{for any $P$}, on the right tail using negative association and the Kearns-Saul inequality. Further, \cite{McAllester2003} showed the left tail bound $\pr(M_0 \leq E[M_0]-\epsilon)\le e^{-1.359n\epsilon^2}$, $\epsilon \ge 0$, using negative association and connections between Chernoff entropy and Gibbs variance. In \cite{berend2013}, the left tail bound was further improved to 
\begin{equation}
    \pr(M_0 \leq E[M_0]-\epsilon) \le e^{-1.92n\epsilon^2},\ \epsilon \ge 0, \label{eq:left_tail_best_sG}
\end{equation}
by showing that $-M_0(X^n, P)\sim\sG(\gamma / n),$ \textcolor{blue}{for any $P$}, on the right tail, where $\gamma$ is defined as
\begin{gather}
    \gamma = \max_{t>0}\,t\,e^{-t}(1-e^{-t}) = 0.2603\ldots,\ (1/2\gamma)= 1.92\ldots.
\end{gather}
\noindent\emph{Remarks}: As we show in Section \ref{subsec: Var_bound}, $\tv(M_{0})\le \gamma/n$, and there exist uniform distributions that nearly achieve this variance upper bound \cite{Rajaraman17}. Now, if $Z\sim\sG(v)$, it is known that $v\ge\text{var}(Z)$ \cite{boucheron2013}. So, the constant 1.92 in the left tail bound $e^{-1.92 n \epsilon^2}$ is nearly optimal. 

In our recent work \cite{Chandra_ncc}, we provided evidence that the constant 1 is likely to be tight for the right tail bound $e^{-n \epsilon^2}$ in the sub-Gaussian regime. More precisely, we showed that, in the Poisson sampling model ($N\sim\text{Poisson}(n)$ samples), there exists $P$ such that $M_0(X^N,P)$ is \textbf{not} $\sG(c/n)$ for $c<0.5$. 

Therefore, the best \textcolor{blue}{distribution free (i.e. uniform over all $P$)} asymptotic left and right tail bounds for $M_0$ using sub-Gaussianity appear to be $e^{-1.92 n \epsilon^2}$ and $e^{-n \epsilon^2}$, respectively. However, it is interesting to consider tail bounds using other approaches such as sub-Poisson and sub-Gamma property \cite{ben-hamou2017}, as these are known to improve upon sub-Gaussian tail bounds in some cases \cite{boucheron2013}. 

Extending concentration results from missing mass $M_0(X^n,P)$ to missing mass of functions $M_{0,g}(X^n,P)$ is challenging because the methods used for proving concentration results for $M_0$ are specific to the missing mass $M_0$ and to sub-Gaussian concentration results. For applications, obtaining the best possible concentration results by going beyond sub-Gaussianity, if necessary, is interesting and important.

\subsubsection{Our Results}\label{sec:conc_results}
To obtain concentration results and tail bounds for $M_{0,g}(X^n,P)$ that are distribution-free, a standard approach is to upper bound the log Moment Generating Function (log MGF) $L_{M_{0,g}}(\lambda)=\log E[e^{\lambda (M_{0,g} - E[M_{0,g}])}]$, by a function of $\lambda$, and use the Chernoff method. Now, $M_{0,g}=\sum_x g(p_x)I(N_x=0)$, is a sum of random variables $g(p_x)I(N_x=0)$, $x\in\cX$, which are (1) dependent because $N_x$ are dependent, and (2) highly heterogeneous because of the range of values of $p_x$. To address dependency, we will follow \cite{McAllester2003} and use the idea of negative association. The log MGF of $M_{0,g}$
\begin{align}
    L_{M_{0,g}}(\lambda) &= \log E\left[\prod_{x\in\cX} \exp\{\lambda g(p_x) 
                        [I(N_x=0)-\pr(N_x=0)]\}\right]\nonumber\\
    &\overset{(a)}{\le} \sum_{x\in\cX} \log E\left[\exp\{\lambda g(p_x) 
                        [I(N_x=0)-\pr(N_x=0)]\}\right],\label{eq:logMGF_na}
\end{align}
where $(a)$ follows because of the following argument. The random variables $N_x$, $x\in\cX$, are negatively associated \cite{McAllester2000, Dubhashi1998} and $$T(N_x)\triangleq\exp\{\lambda g(p_x)[I(N_x=0)-\pr(N_x=0)]\}$$ 
is a monotonically non-increasing function of $N_x$. So, $T(N_x)$, $x\in\cX$, are negatively associated \cite{Dubhashi1998} resulting in $E[T(N_x)T(N_{x'})]\le E[T(N_x)]E[T(N_{x'})]$.

Following the above simplification, to address the heterogeneity and decouple $\lambda$ and $p_x$, \cite{McAllester2000, berend2013} use the Kearns-Saul inequality to upper bound $E[T(N_x)]$. An alternative method, introduced in \cite{ben-hamou2017}, starts with Bennett's inequality \cite{boucheron2013}. We will employ a stronger version of Bennett's inequality for a zero-mean random variable $Z$ with $|Z|<1$, which is proved as follows:
\begin{align}
    E[\exp(tZ)] &\overset{(E[Z]=0)}{=} 1+\sum_{k=2}^{\infty}\frac{t^kE[Z^2Z^{k-2}]}{k!}\nonumber\\
    &\overset{(|Z|<1)}{\le} 1+\sum_{k=2}^{\infty}\frac{t^kE[Z^2]}{k!}\nonumber\\
    &\le 1+\text{Var}(Z)(e^t-t-1).\label{eq:strongBennett}
\end{align}
By using $1+x\le e^x$ on the RHS of \eqref{eq:strongBennett} with $x=\text{Var}(Z)(e^t-t-1)$ and taking logarithms, we obtain the usual Bennett's inequality, given below: 
\begin{equation}
    \log E[\exp(t Z)] \le \text{Var}(Z)(e^t-t-1),\text{ all } t.
    \label{eq:Bennett's}
\end{equation}
Setting $t=\lambda g(p_x)$ and $Z=I(N_x=0)-\pr(N_x=0)$ with $\text{Var}(Z)=(1-p_x)^n(1-(1-p_x)^n)$ in \eqref{eq:strongBennett} and \eqref{eq:Bennett's}, and substituting the upper bounds in \eqref{eq:logMGF_na}, we get the following two bounds on the log MGF:
\begin{align}
    L_{M_{0,g}}(\lambda) &\le \sum_{x\in\cX} \log(1+(1-p_x)^n(1-(1-p_x)^n)(e^{\lambda g(p_x)}-\lambda g(p_x) - 1)), \label{eq:log_MGF_strongbennett}\\
    L_{M_{0,g}}(\lambda) &\le \sum_{x\in\cX} (1-p_x)^n(1-(1-p_x)^n)(e^{\lambda g(p_x)}-\lambda g(p_x) - 1).\label{eq:log_MGF_bennett}
\end{align}
To obtain distribution-free concentration results and tail bounds, the next step is to upper bound the RHS above by a function of $\lambda$ (free of the distribution $p_x$) for $\lambda>0$ to obtain right tail bounds and $\lambda<0$ to obtain left tail bounds. 
\begin{table}[htb]
    \centering
    \caption{Types of concentration and tail bounds.}
    \label{tab:con_tail_bounds}
    \begin{tabular}{|c|c|c|c|c|}
        \hline
        Type & \makecell*[cc]{Upper bounding\\Function} & Range & Tail bound & Remarks \\
        \hline
        General&$f(\lambda)$&$\lambda\in[s,t]$&\makecell*{$\exp\left\{\min\limits_{\lambda\in[s,t]}f(\lambda)-\lambda\epsilon\right\}$}&Chernoff\\
        \hline
        \makecell*[cc]{Sub-Gaussian$(\sigma^2)$\\$\sigma^2$: variance parameter} & $\lambda^2 \sigma^2/2$  & \makecell*[cr]{Right: $\lambda>0$\\ Left: $\lambda<0$} & $e^{-\epsilon^2/2\sigma^2}$ & \cite{boucheron2013}\\
        \hline
        \makecell*[cc]{Sub-Gamma$(v,c)$\\$v$: variance parameter\\$c$: scale parameter} & \makecell*{$\dfrac{\lambda^2v/2}{1-c\lambda}$} & \makecell[cc]{Right:\\$\lambda\in[0,1/c)$} & \eqref{eq:sub-gamma-tail-bound} & \cite{boucheron2013}\\
        \hline
        Strongly-sub-Gamma$(v,c)$ & \makecell*{$\frac{v}{c^2}\log\dfrac{e^{-\lambda c}}{1-c\lambda}$} & \makecell[cc]{Right:\\$\lambda\in[0,1/c)$} & \eqref{eq:strongly-sub-gamma-tail-bound} & Proposed\\
        \hline
        \makecell*[cc]{Filtered sub-Gaussian$(\sigma^2,h(\lambda))$\\$h(\lambda)$: filter} & \makecell*[cc]{$\lambda^2\sigma^2/2+h(\lambda)$}& \makecell[cc]{$\lambda\in[s,t]$} & \makecell*[cc]{$\exp\left\{\min\limits_{\lambda\in[s,t]}\lambda^2\sigma^2/2+h(\lambda)-\lambda\epsilon\right\}$} & Proposed\\
        \hline
        \makecell*[cc]{Poly-sub-Gamma-filtered$(R,\mathbf{a},v,c)$\\$R$: degree parameter\\$\mathbf{a}=[a_2\ a_3\ \cdots\ a_R]$} & \makecell*[cc]{$\sum\limits_{r=2}^{R}a_r\dfrac{\lambda^r}{r}+$\\ $\frac{v}{c^2}\log\dfrac{e^{-\lambda c}}{1-c\lambda}$}& \makecell[cc]{Right:\\$\lambda\in[0,1/c)$} & \makecell*[cc]{$R=2$: \eqref{eq:poly-strongly-sub-gamma-tail-bound}\\$R>2$: numerical} & Proposed\\
        \hline
    \end{tabular}
\end{table}
Table \ref{tab:con_tail_bounds} shows the different forms of the upper-bounding function of $\lambda$ along with the range of $\lambda$ and associated tail bounds obtained by the Chernoff method. The first row describes the general Chernoff bounding method. 

Sub-Gaussianity is a standard notion, and has been extensively studied for the case of missing mass where $g(p)=p$. While sub-Gamma concentration is described in \cite{boucheron2013}, we propose the notion of strongly sub-Gamma concentration, where the log-MGF is upper-bounded by $(v/c^2)\log\frac{e^{-\lambda c}}{1-c\lambda}$, which is the log-MGF of a Gamma random variable with variance parameter $v$ and scale parameter $c$. Since the function $\frac{\lambda^2v/2}{1-c\lambda}$, introduced for sub-Gamma concentration in \cite{boucheron2013}, does not correspond to the log-MGF of a Gamma random variable, the notion of strongly sub-Gamma concentration is, perhaps, a useful addition. It is easy to show the following implications:
\begin{equation}
    Z\sim\sG(v),\,\lambda>0 \to Z\sim\text{Strongly sub-Gamma}(v,c),\,c>0 \to Z\sim\sg(v,c).
\end{equation}
\textcolor{blue}{Note that the variance factor $v$ with which a random variable $Z$ can be $\sG$ or {Strongly sub-Gamma} or $\sg$ is lower bounded by $\text{Var}(Z).$}
So, while sub-Gaussianity is the strongest notion \textcolor{blue}{i.e. it gives the best tail bound for a given variance factor}, the optimality of the variance factor is critical for obtaining the best possible tail bounds. If the variance factor of sub-Gaussian concentration does not coincide with the variance of $Z$, \textcolor{blue}{there is a possibility of} the seemingly weaker notions of concentration, such as strongly sub-Gamma or sub-Gamma concentration \textcolor{blue}{with smaller variance factors than $\sG$}, resulting in better tail bounds on $\pr(Z>E[Z]+\epsilon)$ for a restricted range of $\epsilon$, as illustrated for missing mass later on.

Next, we propose the notion of filtered sub-Gaussian concentration, where the log-MGF is upper-bounded by $f(\lambda)=\lambda^2\sigma^2/2+h(\lambda)$. In terms of distributions, this is equivalent to convolving a Gaussian distribution with another distribution with log-MGF $h(\lambda)$ (assuming it is a valid log-MGF corresponding to some distribution). Hence, we use the terminology of ``filtering'' and call $h(\lambda)$ as a filter. The convolution seemingly weakens concentration when compared to sub-Gaussian concentration. However, when the variance factor of sub-Gaussian concentration is not optimal, instead of bounding the tail with that of a Gaussian random variable of higher variance, one can possibly obtain better bounds by considering convolutions of a Gaussian random variable (with lower variance) and other random variables. 

The last row of the table describes the specific type of filtering that is useful for missing mass, namely, filtering a Gamma log-MGF by a polynomial. When the degree parameter $R=1$, this type of concentration reduces to strongly sub-Gamma concentration. When $R=2$, the log-MGF upper bound $\lambda^2a_2/2+(v/c^2)\log\frac{e^{-\lambda c}}{1-c\lambda}$ corresponds to a random variable whose distribution is the convolution of a Gaussian random variable with variance $a_2$ and a Gamma random variable with variance parameter $v$ and scale parameter $c$. As we illustrate later, such convolutions improve the right tail bound for missing mass. For $R>2$, the polynomial part may result in the filter not being a valid log-MGF, but the tail bounds are still valid, and can be computed numerically.

If $Z\sim\sg(v,c)$ on the right tail, we get the following tail bound using the standard Chernoff method:
\begin{equation}
    \pr(Z\ge E[Z]+\epsilon)\le \exp\left\{-\left(1+\frac{c\epsilon}{v}-\sqrt{1+\frac{2c\epsilon}{v}}\right)\frac{v}{c^2}\right\}.
    \label{eq:sub-gamma-tail-bound}
\end{equation}
Similarly, if $Z\sim\text{Strongly sub-Gamma}(v,c)$, we get
\begin{equation}
    \pr(Z\ge E[Z]+\epsilon)\le 
    \exp\left\{ - \frac{1}{c} \left(\epsilon - \frac{v}{c} \log \left(1 + \frac{c\epsilon} {v}\right) \right)  \right\}.
    \label{eq:strongly-sub-gamma-tail-bound}
\end{equation}
\textcolor{blue}{
Using the Taylor series bound $\ln(1+x) \leq x - 0.5x^2 + (1/3)x^3, x \geq 0,$ on the RHS in \eqref{eq:strongly-sub-gamma-tail-bound}, we get
\begin{equation}
    \pr(Z\ge E[Z]+\epsilon)\le \exp\left\{ -\frac{\epsilon^2}{2v} \left( 1 - ({2}/{3}) \frac{c\epsilon}{v} \right) \right\}  \tag{\ref{eq:strongly-sub-gamma-tail_ub}}
\end{equation}
Note that the exponent in the RHS of \eqref{eq:strongly-sub-gamma-tail_ub} is away by $ (1/3) (c \epsilon^3/v^2)$   from the exponent in the tail bound for $\sG(v).$}
If $Z\sim\text{Poly-sub-Gamma filtered}(R,\mathbf{a},v,c)$, the minimization inside the exponent in Chernoff bound computation is easy to evaluate numerically for any $R$ because $f'(\lambda)$ is rational in $\lambda$ with a linear denominator. For $R=2$, the tail bound can be expressed analytically as follows. 
\begin{align}
&\pr(Z - {E}[Z] \geq \epsilon) \leq \exp\left\{-\frac{1}{c}\left( \left(\frac{1}{2} -\frac{d_2}{d_1}\right)\epsilon +\frac{v}{c} \ln\left( 1-\frac{2c\epsilon}{ d_1 } \right)  \right) \right\}, \label{eq:poly-strongly-sub-gamma-tail-bound}
\end{align}
where $d_1 = (a_2+v+c\epsilon)+\sqrt{(a_2+v+c\epsilon)^2 - 4a_2c\epsilon  }$ and $d_2 =a_2-(v+c\epsilon)$.
\textcolor{blue}{For $a_2 \leq v,$ we upper bound the exponent in the RHS of \eqref{eq:poly-strongly-sub-gamma-tail-bound} with its corresponding third order Taylor series polynomial about $\epsilon = 0,$ to get 
\begin{equation}
    \pr(Z\ge E[Z]+\epsilon)\le \exp\left\{ -\frac{\epsilon^2}{2(v+a_2)} \left( 1 - ({2}/{3}) \frac{cv\epsilon}{(a_2 + v)^2} \right) \right\}  \tag{\ref{eq:poly-strongly-sub-gamma-tail_ub}}
\end{equation} }
Using the definitions and notation in Table \ref{tab:con_tail_bounds}, we now describe our main concentration result for the missing mass of the function $g(p)$, $M_{0,g}(X^n,P)$, in Theorem \ref{thm:con_G0}. Let 
\begin{align}
 u^*_{r}(n,g) &\triangleq \max_{0 < p < 1}  g(p)^r (1-p)^n\, \dfrac{1- (1-p)^n}{p}.\label{eq:ustar}
\end{align}
Where necessary, we will drop the arguments $n,g$ in $u^*_{r}(n,g)$ to reduce clutter.
\begin{theorem}
\begin{enumerate}
\item 
\begin{equation}
    M_{0,g}(X^n,P)\sim \sG\left(\frac{0.519\max_p(g(p)/p)^2}{n}\right),\ \text{all }\lambda \in \mathbb{R}.\label{eq:subGall}
\end{equation}
Additionally, on the left tail, 
\begin{equation}
    M_{0,g}(X^n,P)\sim \sG(u^*_2(n,g)),\ \lambda<0.
\end{equation}
\item Let $R\ge 1$ be a positive integer. Let $c > 0$ be such that 
\begin{equation}
\frac{u^*_{r}(n,g)}{(r-1)!} \leq c \frac{u^*_{r-1}(n,g)}{(r-2)!},\ r\ge3.
\label{eq:scale_parameter}
\end{equation}
Let $v=u^*_{R+1}(n,g)/(c^{R-1}R!)$, and
$$a_r=\left(\frac{u^*_{r}(n,g)}{(r-1)!} - c^{r-2} v\right),\ 2\le r \le R.$$ 
Then, for $\mathbf{a}=[a_2\ a_3\ \cdots\ a_R]$, on the right tail,
\begin{equation}
M_{0,g}(X^n,P) \sim \text{Poly-sub-Gamma-filtered}(R,\mathbf{a},v,c).
\end{equation}
\end{enumerate}
\label{thm:con_G0}
\end{theorem}
The first part of Theorem \ref{thm:con_G0}, proven in Section \ref{subsec:conc_proofs_left}, provides a sub-Gaussian right-tail and left-tail bound as shown in Table \ref{tab:con_tail_bounds}. For the case of missing mass, i.e. $g(p)=p$, we have sub-Gaussianity with a variance of $0.519/n$, which is close to the best-known $0.5/n$. Additionally, since $u^*_2(n,p)\le \gamma/n$ (proved in Section \ref{subsec: Var_bound}; recall that $\gamma = \max_{t>0}\,t\,e^{-t}(1-e^{-t}) = 0.2603\ldots$), this recovers the best-known \textcolor{blue}{distribution free or uniform} sub-Gaussian left-tail bound $e^{-1.92n\epsilon^2}$. Our proof of these tail bounds for $M_0$ are arguably simpler compared to the Gibbs-variance method used in \cite{McAllester2003} or the Kearns-Saul inequality method used in \cite{berend2013}. For all other $g(p)$, such as $g(p)=p^{\alpha}$ for power sum or $g(p)=p\log(1/p)$ for Shannon entropy, the sub-Gaussian concentration left tail results have been shown for the first time. 

For the case of $g(p)=p^{\alpha}$ and $g(p)=p\log_2 (1/p)$, we simplify the left tail bounds and present them in the following corollary.
\begin{corollary}
 \begin{enumerate}
     \item The order-$\alpha$ missing mass $M_{0,\alpha}$ for $\alpha \geq 1,$ satisfies 
     \begin{align}
         &\pr(M_{0,\alpha} - {E}[M_{0,\alpha}] \leq -\epsilon) \leq \exp \left\{ -\frac{n^{2\alpha -1} \epsilon^2}{2 \gamma_\alpha}\right\},\quad n \geq (2\alpha-1) \ln 2/(2\alpha -1-\ln 2),
     \end{align}
     where $\gamma_{\alpha} = \max_{t > 0} t^{2\alpha-1} e^{-t} (1-e^{-t})$.
     \item Let $H_0(X^n, P) \triangleq \sum_{x \in \cX} p_x \log_2 (1/p_x) I(N_x(X^n) =0)$ be the missing Shannon entropy in $X^n$ of a distribution $P.$  $H_0(X^n, P)$ satisfies 
     \begin{align}
         &\pr(H_{0} - {E}[H_{0}] \leq -\epsilon) \leq \exp \left\{ -\frac{n \epsilon^2}{2 \gamma (\log_2 n)^2}  \right\},\quad n \geq 3,
     \end{align}
     where $\gamma=\max_{t>0} t e^{-t}( 1-e^{-t}) = 0.2603\ldots$.
 \end{enumerate}
 \label{cor:left_tail_special_case}
\end{corollary}
\textcolor{blue}{The variance factors in the above uniform sub-Gaussian left tail bounds for $M_{0,\alpha}$ and $H_0$ are near optimal.}
A proof of Corollary \ref{cor:left_tail_special_case} is given in Section \ref{sec:corollary_proof}. 

The second part of Theorem \ref{thm:con_G0}, proven in Section \ref{subsec:conc_proofs_right}, provides right tail bounds using the Poly-sub-Gamma-filtered approach with a Gamma MGF and a polynomial  filter as described in Table \ref{tab:con_tail_bounds}. The choice of the scale parameter $c$ needs elaboration. From \eqref{eq:scale_parameter}, it is clear that $c$ quantifies the rate of decay of $u^*_r(n,g)/(r-1)!$, and it will depend on the function $g(p)$. An important requirement is that $c$ should decay with $n$ so that the confidence intervals from the tail bounds shrink with increasing $n$. The following \emph{rate of decay} lemma describes the choice of $c$ as a decreasing function of $n$ for two types of functions $g(p)$.

\begin{lemma}
Let $g(p)$ be differentiable for $p\in(0,1)$, and let $g'(p)$ denote its derivative. Let $u^*_r(n,g)$ be as defined in \eqref{eq:ustar}. There exist $b_i$, $0\le i\le 5$, independent of $r$, such that
\begin{equation}
    c= \max( b_0 \ g(b_1/(n+b_2)), b_3 \ g(b_4/(n+b_5)))\text{ satisfies }\frac{u^*_{r}(n,g)}{(r-1)!} \leq c \frac{u^*_{r-1}(n,g)}{(r-2)!},\ r\ge3, \label{eq:rate_of_fall} \end{equation}
if $g(p)$ belongs to one of the following two types:
\begin{enumerate}
    \item Type A: For some $\mu>0$, 
    $$0 < g'(p) \ \le \ \mu \ g(p)/p,\ p\in(0,1).$$
    \item Type B: For some $p^*\in(0,1)$, 
    \begin{align*}
        0 < g'(p) &\le (1/p - 1/p^*)\dfrac{g(p)}{1-p},\ p\in(0,p^*),\\
        g'(p) &< 0,\ \phantom{(1p - 1/p^*)\dfrac{g(p)}{1-p}} p\in(p^*,1).
    \end{align*}
\end{enumerate}
\label{lem:rec}
\end{lemma}
A proof of Lemma \ref{lem:rec} is given in Section \ref{sec:rec_proof}. In the proof, the values of $b_i$ are explicitly specified for a $g(p)$ of Type A or B. Interesting examples such as  $g(p)=p^\alpha$, $\alpha>0$, belong to Type A with $\mu=\alpha$ and $g(p)=p\log_2 (1/p)$ belongs to Type B with $p^*=1/e$. For these cases, we further simplify the right tail bound for $R=1$ (strongly sub-gamma concentration) in the following corollary.
\begin{corollary}
 \begin{enumerate}
     \item The missing mass $M_0$ satisfies
\begin{align}
\pr(M_{0} - {E}[M_{0}] \geq \epsilon) & \overset{(a_1)}{\leq} \exp \left\{ - \frac{2(n+2)}{3} \left( \epsilon - \frac{2\gamma_n} {3} \log \left(1 + \frac{   3 \epsilon }{ 2\gamma_n } \right)  \right) \right\},  n \geq 3, \label{eq:right_tail_derived_ssg_M0} \\
& \textcolor{blue}{\overset{(a_2)}{\leq} \exp \left\{ -1.92n\epsilon^2 (1 - (\epsilon /\gamma_n) ) \right\}, n \geq 3,} \label{eq:right_tail_derived_ssg_M0_ub} \\ 
\pr(M_{0} - {E}[M_{0}] \geq \epsilon) & \textcolor{blue}{\overset{(a_3)}{\leq} \exp \left\{ -1.92n\epsilon^2 (1 - ( \psi/3\gamma^2) \epsilon ) \right\}, n \geq 3,} \label{eq:right_tail_derived_poly_ssg_M0_ub}
\end{align}
where $\gamma_n = \gamma (1+2/n)$, $\gamma=\max_{t>0} t e^{-t}( 1-e^{-t}) = 0.2603\ldots,$ and $\psi = \max_{t>0} t^2 e^{-t}( 1-e^{-t}) = 0.477\ldots $. 
\label{corr_M0: ssg_M0}
\item Consider the order-$\alpha$ missing mass $M_{0,\alpha}$, $\alpha > 1$. For $n>1+4\alpha^2/(1-\alpha)^2$, 
\begin{align}
    &\pr(M_{0,\alpha} - {E}[M_{0,\alpha}] \geq \epsilon) \leq \exp \left\{-\frac{n-b}{a}\epsilon +\frac{\gamma_{\alpha}}{a^2} (n-b)^{3-2\alpha} \log \left(1+\frac{a n^{2\alpha -1}}{\gamma_{\alpha} (n-b)}\epsilon\right) \right\}, \label{eq:right_tail_derived_ssg_M0_alpha}
\end{align}
where $b = 1+\frac{2\alpha}{\alpha -1}$, $a = (b-1)\left( \frac{2(\alpha-1)}{\alpha +1} \right)^{\alpha}$ and $\gamma_{\alpha} = \max_{t > 0} t^{2\alpha-1} e^{-t} (1-e^{-t})$. 
\item For the missing Shannon entropy in $X^n$, $H_0(X^n, P)$, of a distribution $P,$
\begin{align}
    \pr(H_{0} - {E}[H_{0}] \geq \epsilon) & \overset{(b_1)}{\leq}  \ \exp \bigg\{ -\frac{2n_0}{\log_2 n_0} 
    \bigg(\epsilon -  \frac{ \gamma'_{n} }{ \log_2(n_0) } \log \left( 1 + \epsilon \frac{\log_2(n_0)}{ \gamma'_{n}   } \right)   \bigg) \bigg\}, n \geq 3, \label{eq:right_tail_derived_ssg_H_0} \\
    & \textcolor{blue}{\overset{(b_2)}{\leq} \ \exp \bigg\{ -\frac{n_0 \epsilon^2}{\gamma'_n} 
    \bigg(1 - (2/3) \frac{ \epsilon \log_2(n_0) }{ \gamma'_n } \bigg) \bigg\}, n \geq 3,} \label{eq:right_tail_derived_ssg_H_0_ub} 
\end{align}
where  $n_0 = (1/3)(n-1) + e$ and  $\gamma'_{n} = 2 \gamma (1/3 + (e-1/3)/n)   (\log_2 n)^2$.
 \end{enumerate}
 \label{cor:right_tail_special_case}
\end{corollary}
 \textcolor{blue}{Note that we get \eqref{eq:right_tail_derived_ssg_M0} by showing that $M_0$ is strongly sub-gamma on the right tail with variance factor 
 $\gamma/n$ and scale parameter $3/2(n+2)$ i.e. we have shown that $M_0$ is uniformly strongly sub-gamma on the right tail with a near-optimal variance factor. We get the bound in \eqref{eq:right_tail_derived_ssg_M0_ub} by using the right tail bound in \eqref{eq:strongly-sub-gamma-tail_ub} for a strongly sub-gamma$(\gamma/n, 3/2(n+2))$ random variable.  We get the bound in  \eqref{eq:right_tail_derived_poly_ssg_M0_ub} by showing that $M_0$ is Poly-sub-Gamma-filtered with $R=2, v = (1+2/n) \ \psi/3n, c = 3/2(n+2),$ and $ a = [a_2] =  [\gamma/n - v].$
 The bounds in \eqref{eq:right_tail_derived_ssg_M0_ub}, \eqref{eq:right_tail_derived_poly_ssg_M0_ub} are close to $e^{-1.92n\epsilon^2}$ (the best sub-gaussian left tail bound for $M_0$) for small $\epsilon.$ 
 The variance factors in the uniform strongly sub-gamma results for the order-$\alpha$ missing mass and the missing Shannon entropy in \eqref{eq:right_tail_derived_ssg_M0_alpha}, \eqref{eq:right_tail_derived_ssg_H_0} respectively are also nearly-optimal. The bounds obtained for $R=2,3\ldots,$ also have near-optimal variance factors that are non-increasing.} A proof of Corollary \ref{cor:right_tail_special_case} is given in Section \ref{sec:corollary_proof}.
    
\begin{figure}[tbh]
\begin{center}
\input{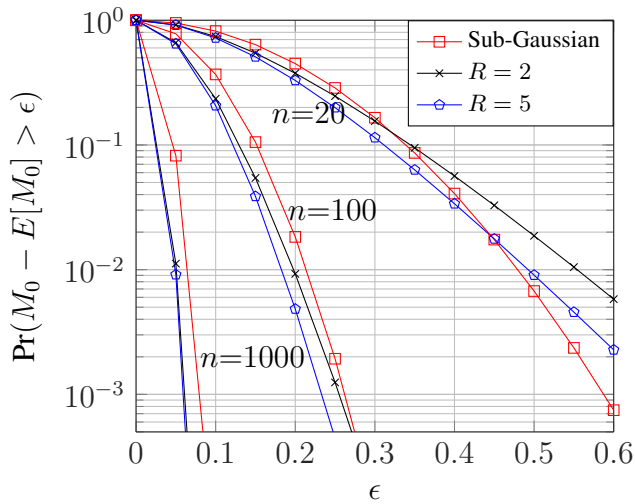}
\end{center}
  \caption{Comparison of right tail bounds for $M_0$.}
  \label{fig:plot1}
\end{figure}

For $M_0$ (i.e. $g(p) = p$), Fig. \ref{fig:plot1} shows a plot comparing the different right tail bounds for $n=20,100,1000$.
The new bounds, for $R=2$ and $R=5$, provide a noticeable improvement over the sub-Gaussian bound $e^{-n\epsilon^2}$ for all three values of $n$ over a significant range of $\epsilon$. 


The proofs of the results in this section can be found in Section \ref{sec:minimax_M_0_alpha_proof} and in the ensuing ones.

\section{Proof of Theorem \ref{thm:minimax_M_0_alpha}} \label{sec:minimax_M_0_alpha_proof}
We start by proving the upper bound of Part 1 of the theorem.
\subsection{Upper bound on $R^*_{n,\alpha}$ for $\alpha \in \mathbb{N}$}
Consider the generalized Good-Turing estimator $\mgt_{0,\alpha}(X^n) = \frac{\phi_{\alpha}(X^n)}{\binom{n}{\alpha}}$, $\alpha \in \mathbb{N}$. By the definition of minimax risk, we have  
\begin{align}
R_{n,\alpha}^* & = \min_{\widehat{M}_{0,\alpha}} \max_{P} E [({\widehat{M}_{0,\alpha}}(X^n) - M_{0,\alpha}(X^n,P))^2] \nonumber \\ & \leq  \max_{P} E [(\mgt_{0,\alpha}(X^n) - M_{0,\alpha}(X^n,P))^2]\triangleq R_{n,\alpha}(\mgt_{0,\alpha}). \label{eq:ub_method}
\end{align}
The rest of the proof is to upper bound the worst-case squared-error risk of the generalized Good-Turing estimator $R_{n,\alpha}(\mgt_{0,\alpha})$, and this is provided in the next lemma.
\begin{lemma}
 \begin{align}
       &R_{n,\alpha}(\mgt_{0,\alpha}) \leq \ O\left(1/n^{2\alpha -1}\right) . \label{eq:worst_case_G0_alpha}
 \end{align}
 \label{lem:minimax_ub_int_M_0_alpha}
\end{lemma}
\begin{proof}
\begin{align}
    (\mgt_{0,\alpha} - M_{0,\alpha})^2 &= 
    \left(\sum_{x \in \mathcal{X}} \frac{1}{\binom{n}{\alpha}} I(N_x = \alpha) - p_x^{\alpha}I(N_x = 0) \right)\left(\sum_{y \in \mathcal{X}} \frac{1}{\binom{n}{\alpha}} I(N_y = \alpha) - p_y^{\alpha}I(N_y = 0) \right) \nonumber \\ 
    & = \Bigg[\sum_{x \in \mathcal{X}} \bigg( \frac{1}{{\binom{n}{\alpha}}^2} I(N_x = \alpha) + p_x^{2\alpha}I(N_x = 0) \bigg)\Bigg] \nonumber \\ 
    & \quad + \Bigg[\sum_{x \in \mathcal{X}} \sum_{y \in \mathcal{X}, y \neq x} \frac{1}{{\binom{n}{\alpha}}^2} I(N_x = N_y =  \alpha) - \sum_{x \in \mathcal{X}} \sum_{y \in \mathcal{X}, y \neq x} \frac{p_x^{\alpha}}{\binom{n}{\alpha}} I(N_x =0, N_y = \alpha)\Bigg] \nonumber \\ 
    & \quad - \Bigg[\sum_{x \in \mathcal{X}} \sum_{y \in \mathcal{X}, y \neq x} \frac{p_y^{\alpha}}{\binom{n}{\alpha}} I(N_y =0, N_x = \alpha) - \sum_{x \in \mathcal{X}} \sum_{y \in \mathcal{X}, y \neq x} (p_x p_y)^{\alpha} I(N_x = N_y =  0)\Bigg] \label{eq:square_brackets}
    \end{align}
Taking expectations, we consider each of the terms inside square brackets above and bound them by $O(1/n^{2\alpha-1})$. 

The following bounds are used repeatedly in the simplifications. Let $f_1:[0,1]\to\mathbb{R}^+$ and $f_2:[0,1]^2\to\mathbb{R}^+$ be non-negative functions, and let $a$, $b$, $c$ be positive numbers.
\begin{align}
\sum_x p_x f_1(p_x) &\le \max_{p\in [0,1]} f_1(p) \label{eq:avg<max1}\\
\sum_x\sum_y p_x p_y f_2(p_x,p_y) &\le \max_{p,q\in[0,1]} f_2(p,q) \label{eq:avg<max2}\\
\max_{p\in[0,1]} p^a(1-p)^b &= \dfrac{a^a b^b}{(a+b)^{a+b}} \label{eq:max2}\\
\max_{\substack{p,q\in[0,1]\\p+q\le1}} p^a q^b (1-p-q)^c &= \dfrac{a^a b^b c^c}{(a+b+c)^{a+b+c}} \label{eq:max3}
\end{align}
Let us first consider the first term in square brackets in \eqref{eq:square_brackets}.
\begin{align} 
    \sum_{x \in \mathcal{X}}  E\left( \frac{1}{{\binom{n}{\alpha}}^2} I(N_x = \alpha) + p_x^{2\alpha}I(N_x = 0) \right)  
    &= \sum_{x \in \mathcal{X}} \frac{1}{{\binom{n}{\alpha}}}  p_x^{\alpha} (1-p_x)^{n-\alpha} + \sum_{x \in \mathcal{X}} p_x^{2\alpha} (1-p_x)^{n} \nonumber \\ 
    & \overset{(a)}{\leq} \frac{1}{{\binom{n}{\alpha}}} \frac{(\alpha-1)^{\alpha-1} (n-\alpha)^{n-\alpha}}{(n-1)^{n-1} } + \frac{(2\alpha-1)^{2\alpha-1} (n)^{n}}{(n+2\alpha-1)^{n+2\alpha-1} } \nonumber \\ 
    & \leq O\left(\frac{1}{n^{2\alpha-1}}\right), \label{eq:b11}
   \end{align} 
where $(a)$ uses \eqref{eq:avg<max1}, \eqref{eq:max2}. 
Next, we consider the second term in square brackets in \eqref{eq:square_brackets}. We will need the notation $\binom{n}{a,b}=\frac{n!}{a!b!(n-a-b)!}$. Further, we will assume $n>2\alpha$.
\begin{align} 
    &\sum_{x \in \mathcal{X}}\sum_{y \in \mathcal{X}, y \neq x} E \bigg( \frac{1}{{\binom{n}{\alpha}}^2} I(N_x = N_y =  \alpha) - \frac{p_x^{\alpha}}{\binom{n}{\alpha}} I(N_x =0, N_y = \alpha)\bigg)\nonumber\\ 
    &\qquad\qquad\qquad\qquad= \sum_{x \in \mathcal{X}}\sum_{y \in \mathcal{X}, y \neq x} (p_x p_y)^{\alpha} \bigg(\frac{\binom{n}{\alpha, \alpha}}{{\binom{n}{\alpha}}^2}  (1-p_x-p_y)^{n-2\alpha} -  (1-p_x-p_y)^{n-\alpha}\bigg) \nonumber \\  
    &\qquad\qquad\qquad\qquad\overset{(a)}{\leq} \sum_{x \in \mathcal{X}}\sum_{y \in \mathcal{X}, y \neq x} (p_x p_y)^{\alpha} (1-p_x-p_y)^{n-2\alpha}(1- (1-p_x-p_y)^{\alpha}) \nonumber \\ 
    &\qquad\qquad\qquad\qquad\overset{(b)}{=} \sum_{x \in \mathcal{X}}\sum_{y \in \mathcal{X}, y \neq x} (p_x p_y)^{\alpha} (p_x+p_y) \left(\sum_{i =0}^{\alpha -1} (1-p_x-p_y)^{n-2\alpha+i}\right) \nonumber \\ 
    &\qquad\qquad\qquad\qquad\overset{(c)}{=} 2\sum_{i=0}^{\alpha -1}\sum_{x \in \mathcal{X}}\sum_{y \in \mathcal{X}, y \neq x} p_xp_y \Big(p_x^\alpha p_y^{\alpha-1}(1-p_x-p_y)^{n-2\alpha+i}\Big) \nonumber \\ 
    &\qquad\qquad\qquad\qquad\overset{(d)}{\leq} 2 \sum_{i =0}^{\alpha -1} \frac{\alpha^{\alpha}(\alpha-1)^{\alpha-1}(n-2\alpha+i)^{n-2\alpha+i}}{(n+i-1)^{n+i-1}} \nonumber \\
    &\qquad\qquad\qquad\qquad \leq O\left(\frac{1}{n^{2\alpha-1}}\right),  \label{eq:b21}
\end{align}
where $(a)$ follows because $\binom{n}{\alpha,\alpha}\le \binom{n}{\alpha}^2$, $(b)$ follows using the identity $1-(1-a)^m=a\sum_{i=0}^{m-1}(1-a)^i$, $(c)$ follows because of the symmetry between $x$ and $y$ in the summation, and $(d)$ uses \eqref{eq:avg<max2}, \eqref{eq:max3}.

Finally, we consider the third term in square brackets in \eqref{eq:square_brackets}.
\begin{align}
    &\sum_{x \in \mathcal{X}}\sum_{y \in \mathcal{X}, y \neq x} E \bigg( \frac{p_y^{\alpha}}{\binom{n}{\alpha}} I(N_x =\alpha, N_y = 0) - (p_x p_y)^{\alpha} I(N_x = N_y =  0)   \bigg) \nonumber \\ 
    &\qquad\qquad\qquad\qquad = \sum_{x \in \mathcal{X}}\sum_{y \in \mathcal{X}, y \neq x} (p_x p_y)^{\alpha} (1-p_x-p_y)^{n-\alpha}(1-(1-p_x-p_y)^{\alpha}) \nonumber \\ 
    &\qquad\qquad\qquad\qquad\leq O\left(\frac{1}{n^{2\alpha-1}}\right), \label{eq:b31}
\end{align}
where the final inequality follows using steps similar to those used to obtain \eqref{eq:b21}. Using \eqref{eq:b11}, \eqref{eq:b21} and \eqref{eq:b31}, the proof is complete.
\end{proof}

\subsection{Lower bound on minimax risk $R^*_{n,\alpha}$}
The lower bound of $O\left({1}/{n^{2\alpha -1}}\right)$ for all $\alpha \in \mathbb{N}$ is obtained by following the Dirichlet prior approach. This approach was used in \cite{Rajaraman17} for the case of missing mass $M_0$ with $\alpha=1$. The same approach is generalized here for $\alpha \in \mathbb{N}$.  

Let $\Delta_k$ be the set of all probability distributions on the alphabet $\mathcal{X} = \{1,2,\ldots,k\}$. Let $P$ be a random variable on $\Delta_k$, generated according to a Dirichlet distribution \cite{ng2011dirichlet} with parameters $\bbeta = (\beta_1, \beta_2, \ldots, \beta_k)$. From the family of distributions $\Delta_k$, $P$ is chosen according to Dirichlet$(\bbeta)$ and then $X^n\sim P$ is sampled \emph{iid}. By the standard prior method, the minimax rate $R^*_{n,\alpha}$ is lower bounded as 
\begin{align} 
R^*_{n,\alpha} &\ge E_{X^n} \left[\text{Var}[M_{0,\alpha}(X^n,P)\vert X^{n}]\right], \label{eq:prior_method}
\end{align}
where $\text{Var}\left[M_{0,\alpha}(X^{n},P)\vert X^{n}\right]$ is the conditional variance of $M_{0,\alpha}(X^{n},P)$ conditioned on $X^n$.
\begin{lemma}
Let $P\sim\text{Dirichlet}(\beta_1,\beta_2,\ldots,\beta_k)$ with $k=cn^2$, $c>0$, and $\beta_i=1/n$. Let $X^n\sim P$, \emph{iid}. Then,
\begin{equation}
E_{X^n} \left[\text{Var} \left[M_{0,\alpha}(X^{n},\ P)\vert X^{n}\right] \ \right]   \ \geq_n \ O\left(1/n^{2\alpha -1}\right),
   \label{eq:minimax_lower_bound}
\end{equation} 
where the notation $\geq_n$ is as defined in Section \ref{sec:prior}.
\label{lem:DirichletPrior}
\end{lemma}
\begin{proof}
A proof is given in Appendix \ref{appsec:dirichletproof}.
\end{proof}
Combining Lemma \ref{lem:DirichletPrior} and \eqref{eq:prior_method}, the proof of the lower bound is complete.

\section{Proof of Theorem \ref{thm:con_G0}}
\label{sec:conc_proofs}
\subsection{Proof of Part 1 of Theorem \ref{thm:con_G0}} 
\label{subsec:conc_proofs_left}
Starting with \eqref{eq:log_MGF_strongbennett}, we simplify as follows:
\begin{align}
    L_{M_{0,g}}(\lambda)& \overset{(a)}{\le} \sum_x \log(1+e^{-np_x}(e^{\lambda g(p_x)}-\lambda g(p_x) - 1)) \nonumber\\
    &= \frac{\lambda^2}{n}\sum_x p_x\left[\frac{n}{\lambda^2p_x}\log(1+e^{-np_x}(e^{\lambda g(p_x)}-\lambda g(p_x) - 1))\right]\nonumber\\
    &\overset{(b)}{\le} \frac{\lambda^2}{n}\max_p\left[\frac{n}{\lambda^2p}\log(1+e^{-np}(e^{\lambda g(p)}-\lambda g(p) - 1))\right]\nonumber\\
    &= \frac{\lambda^2}{n}\max_p\left[\left(\frac{g(p)}{p}\right)^2\frac{np}{(\lambda g(p))^2}\log(1+e^{-np}(e^{\lambda g(p)}-\lambda g(p) - 1))\right]\nonumber\\
    &\le \frac{\lambda^2}{n}\max_p\left(\frac{g(p)}{p}\right)^2\ \max_{u,v}\frac{u}{v^2}\log(1+e^{-u}(e^v-v-1))\nonumber\\
    &\overset{(c)}{\le} \frac{0.2595\lambda^2\max_p(g(p)/p)^2}{n},\label{eq:logMGFsimplified}
\end{align}
where $(a)$ follows by using $(1-p_x)^n(1-(1-p_x)^n)\le (1-p_x)^n \le e^{-np_x}$, $(b)$ follows because average is lesser than maximum, and $(c)$ follows because $\max_{u,v}\frac{u}{v^2}\log(1+e^{-u}(e^v-v-1))=0.2595...$ using some calculus and computations. From \eqref{eq:logMGFsimplified},  we get the sub-Gaussianity result in \eqref{eq:subGall}.

Recall, from \eqref{eq:log_MGF_bennett}, the upper bound
 $$L_{M_{0,g}}(\lambda) \le \sum_{x\in\cX} (1-p_x)^n(1-(1-p_x)^n)(e^{\lambda g(p_x)}-\lambda g(p_x) - 1).$$
Since we are concerned with the left tail bound, we consider $\lambda<0$. 
Using $e^t - t -1 \leq t^2/2$, $t \leq 0$, and setting $t=\lambda g(p_x)$, we get 
\begin{align}
   L_{M_{0,g}}(\lambda) &\leq  \sum_{x \in \cX} (1-p_x)^n (1- (1-p_x)^n) \frac{(\lambda g(p_x))^2}{2} \nonumber \\ 
   &\leq \frac{\lambda^2 u^*_{2}(n,g) }{2},\quad \lambda\le 0, \label{eq:left_tail_mgf_ub}
\end{align}
where $u^*_{2}(n,g)$ is as defined in \eqref{eq:ustar}. This completes the proof of Part 1 of Theorem \ref{thm:con_G0}.
 
 \subsection{Proof of Part 2 of Theorem \ref{thm:con_G0}}
 \label{subsec:conc_proofs_right}
For right tail bounds, we consider $\lambda>0$. The following lemma provides a power series upper bound on $L_{M_{0,g}}(\lambda)$ and is critical in the proof of Part 2 of Theorem \ref{thm:con_G0}.
 \begin{lemma}[Power series upper bound] 
 \begin{equation}
 L_{M_{0,g}}(\lambda) \leq \sum_{r=2}^{\infty} \frac{\lambda^r}{r!} u_r^*(n,g),\quad \lambda>0,
     \label{eq:power_series_ub}
     \end{equation} 
where $ u^*_{r}(n,g)$ is as defined in \eqref{eq:ustar}.
 \label{lem:power_series_ub}
 \end{lemma}
 \begin{proof}
 Using $e^{t} - t - 1 = \sum_{r=2}^{\infty} \frac{t^r}{r!}, t \geq 0,$ with $t = \lambda g(p_x)$ in \eqref{eq:log_MGF_bennett}, we get 
 \begin{align}
      L_{M_{0,g}}(\lambda) &\leq  \sum_{x \in \cX} (1-p_x)^n (1- (1-p_x)^n) \left( \sum_{r=2}^{\infty} \frac{(\lambda g(p_x))^r}{r!} \right) \nonumber \\
      & = \sum_{r=2}^{\infty} \frac{\lambda^r}{r!} \left( \sum_{x \in \cX} g(p_x)^r (1-p_x)^n (1- (1-p_x)^n) \right) \nonumber \\  & = \sum_{r=2}^{\infty} \frac{\lambda^r}{r!} \left( \sum_{x \in \cX} p_x \  g(p_x)^r (1-p_x)^n \frac{(1- (1-p_x)^n)}{p_x} \right) \nonumber \\ & \overset{(a)}{\leq} \sum_{r=2}^{\infty} \frac{\lambda^r}{r!} u_r^*(n,g),\quad \lambda \geq 0, \nonumber
 \end{align}
 where (a) follows by using \eqref{eq:avg<max1} for $\lambda \geq 0$. 
 \end{proof}
Let $R \geq 1.$ For $r\ge R+2$, repeatedly applying $\frac{u^*_{r}(n,g)}{(r-1)!} \le c\, \frac{u^*_{r-1}(n,g)}{(r-2)!},\,r\ge 3 $, we have
\begin{equation}
    \frac{u^*_{r}(n,g)}{(r-1)!} \le c^{r-(R+1)}\, \frac{u^*_{R+1}(n,g)}{R!},\,r\ge R+1.
\end{equation}
Omitting the arguments $n,g$ from $u^*_{r}(n,g)$ for brevity and using the above inequality in \eqref{eq:power_series_ub}, we get
\begin{align}
L_{M_{0,g}}(\lambda) & \leq \sum_{r=2}^{R} \frac{\lambda^r}{r!}\,u^*_{r} + \sum_{r=R+1}^{\infty} \frac{\lambda^r}{r}\,c^{r-(R+1)}\, \frac{u^*_{R+1}}{(R)!},\nonumber\\
&= \sum_{r=2}^{R} \frac{\lambda^r}{r!}\,u^*_{r} + \frac{u^*_{R+1}}{c^{R+1} R!} \sum_{r=R+1}^{\infty}\frac{(c\lambda)^r}{r},\nonumber\\
&\overset{(a)}{=}\left(\sum_{r=2}^{R} \frac{\lambda^r}{r!}\,(u^*_{r}-c^{r-2} v (r-1)!)\right) - \frac{v}{c^2}(c\lambda) + \frac{v}{c^2} \left(\sum_{r=1}^{\infty}\frac{(c\lambda)^r}{r}  \right),    
\end{align}
where $(a)$ follows by setting $v\triangleq\frac{u^*_{R+1}}{c^{R-1} R!}$ and adding/subtracting terms between the two summations. Using the Taylor expansion $-\log(1-x)=\sum_{r=1}^{\infty}\frac{x^r}{r}$, $x<1$, completes the proof.

\section{Proof of Lemma \ref{lem:rec}} 
\label{sec:rec_proof} 
In this section, we present a proof for Lemma \ref{lem:rec}. Let 
$$u_r(p,n,g) \triangleq g(p)^r (1-p)^n (1-(1-p)^n)/p$$ 
be the function that is maximized over $p \in (0,1)$ in \eqref{eq:ustar}, and let $p^*_r=\argmax_{p \in (0,1)}u_r(p,n,g)$. So, we have $u_r^*(n,g) = u_r(p^*_r,n,g)$. In the next lemma, we upper bound the rate of fall of $u_r^*(n,g)/(r-1)!$ using $p^*_r$.
\begin{lemma}
\begin{equation}
    \frac{u_r^*(n,g)}{(r-1)!} \ \leq \ \frac{g(p^*_r)}{r-1} \  \frac{u_{r-1}^*(n,g)}{(r-2)!}, r \geq 3.  \label{eq:c_ub_r}
\end{equation}
\end{lemma}
\begin{proof}
By the definition of $u_{r}(p,n,g), u_{r}^*(n,g)$ and $p^*_{r}$ we get
\begin{align}
    \frac{u^*_{r}(n,g)}{(r-1)!} \ & = \  \frac{u_{r}(p^*_{r},n,g)}{(r-1)!} \ \overset{(a)}{=} \  \frac{g(p^*_{r})}{r-1} \  \frac{u_{r-1}(p^*_{r},n,g)}{(r-2)!} \  \overset{(b)}{\le} \  \frac{g(p^*_{r})}{r-1} \ \frac{u^*_{r-1}(n,g)}{(r-2)!}, \nonumber
\end{align}
where $(a)$ follows from $u_{r}(p,n,g) = g(p) u_{r-1}(p,n,g)$ and $(b)$ follows from $ u_{r-1}(p^*_{r},n,g) \leq u^*_{r-1}(n,g)$.
\end{proof}
In the rest of the proof, we find an upper bound on $c(r) \triangleq \frac{g(p^*_{r})}{r-1}$ that is independent of $r$ and decreasing with $n$ for functions $g(p)$ that are differentiable for $p \in (0,1)$ and fall under either Type A or Type B as described in Lemma \ref{lem:rec}. Recall that $g'(p)$ denotes the derivative of $g(p)$.


\noindent\textbf{Type A}: For this type, we have that for some $\mu>0$, 
    $$0 < g'(p) \le \mu \  g(p)/p,\ p\in(0,1).$$
Taking partial derivative of $u_{r}(p,n,g)$ with respect to $p$, letting $q\triangleq 1-p$ and simplifying, we get
    \begin{align}
        \pdv{ u_{r}(p,n,g)}{p} &= \frac{q^{n-1}}{p^2}\,g(p)^{r-1} \bigg( n pq^n g(p) + (1-q^n)  \bigg[((r-1)q- np) g(p) + rpq \left(g'(p) -\frac{g(p)}{p} \right) \bigg] \bigg). \label{eq:partialderivative} 
     \end{align}
Using   $g'(p)\leq \mu\frac{g(p)}{p}$ in the above equation, we get  \begin{align}  
\pdv{u_{r}(p,n,g)}{p} &\leq \frac{q^{n-1}}{p^2}\,g(p)^r \big( npq^n + (1-q^n)   \big[r\mu-1-(n+r\mu-1)p\big]\big) \nonumber \\ 
& \overset{(a)}{<} 0, \text{ if } r\mu-1-(n+r\mu-1)p < -1, 
\label{eq:u_r_der_A}  
\end{align}
where $(a)$ follows because $npq^n-(1-q^n)<0$. So, $u_r(p,n,g)$ decreases for $p>r\mu/(n+r\mu-1)$ implying that
$$p^*_r\leq \frac{r\mu}{n+r\mu-1}.$$
Since $g(p)$ is an increasing function ($g'(p)>0$),
    \begin{equation}
        c(r)=\frac{g(p^*_r)}{r-1} \leq \frac{1}{r-1} \   g\left(\frac{r \mu  }{n+r\mu-1}\right), r \geq 3. \label{eq:c_ub_r_A}
    \end{equation}
In the next lemma, we present upper bounds on $\frac{1}{r-1} \ g\left(\frac{r \mu  }{n+r\mu-1}\right)$ that are independent of $r$ and are in the form of \eqref{eq:rate_of_fall} for different ranges of $\mu$ and $n$.

\begin{lemma} \label{lem:type-A-c}
     \begin{enumerate}
         \item If either $0 < \mu \leq 1$, or $\mu >1$ and $n < 1+\frac{4\mu^2}{(\mu-1)^2}$,
         $$c(r) \leq 0.5 \ g\left( \frac{3\mu}{n+3\mu-1}\right).$$ 
         To get the form of \eqref{eq:rate_of_fall}, set $b_0 = b_3 = 0.5$, $b_1 = b_4 = 3\mu$ and $ b_2 = b_5 = 3\mu-1$. \label{lem-part:mu < 1}

         \item For $\mu >1$ and $n \geq 1+\frac{4\mu^2}{(\mu-1)^2}$, 
         $$c(r) \leq \max \left( 0.5\,g\left( \frac{3\mu}{n+3\mu-1}\right), \frac{1}{r_2 -1 }\,g\left(\frac{r_2\mu}{n+r_2\mu-1} \right)\right),$$
         where $r_2 = 0.5(n-1)(1-\frac{1}{\mu})\left(1 + \sqrt{1 - \frac{4\mu^2}{(n-1) (\mu - 1)^2}} \right)$. The above bound is directly in the form of \eqref{eq:rate_of_fall}. \label{lem-part:mu>1_disc>0}
     \end{enumerate}
\end{lemma}
\begin{proof}
A proof is given in Appendix \ref{subsec:type-A-c_proof}.
\end{proof}
     
\noindent\textbf{Type B}: For this type, we have that, for some $p^*\in(0,1)$, \begin{align*}
        0 < g'(p) &\le (1/p - 1/p^*)\dfrac{g(p)}{1-p},\ p\in(0,p^*),\\
        g'(p) &< 0,\ \phantom{(1p - 1/p^*)\dfrac{g(p)}{1-p}} p\in(p^*,1).
\end{align*}
We rewrite the partial derivative of $u_{r}(p,n,g)$ in \eqref{eq:partialderivative} as
\begin{align}
    \pdv{u_{r}(p,n,g)}{p} &= \frac{q^{n-1}}{p^2}\,g(p)^{r-1}\bigg(np\,g(p) (2q^n - 1) + (1-q^n)q\bigg[(r-1)g(p)+ rp\left( g'(p) -\frac{g(p)}{p} \right)\bigg]  \bigg). \label{eq:u_der2}  
       \end{align}
We will consider the two ranges $p\in(0,p^*]$ and $p\in(p^*,1)$ separately. For $p\in(p^*,1)$, using $g'(p) < 0$ in \eqref{eq:u_der2} and simplifying, we get
\begin{align}
\pdv{u_{r}(p,n,g)}{p} &\leq \frac{q^{n-1}}{p^2}\,g(p)^r \big(np(2q^n-1) - (1-q^n)q \big)<0, 
\label{eq:u_r_der_B_1}
\end{align}  
where the last inequality can be readily verified. So, $u_r(p,n,g)$ decreases for $p>p^*$ implying that $p^*_r \leq p^*$.

For $p \in (0,p^*]$, using $g'(p) \leq \frac{g(p)}{1-p}\left(\frac{1}{p} - \frac{1}{p^*}\right)$ in \eqref{eq:u_der2} and simplifying, we get 
\begin{align}  
\pdv{u_{r}(p,n,g)}{p} 
&\leq \frac{q^{n-1}}{p^2}\,g(p)^r 
\big(\big[r-1- p(n+(r/p^*)-1)\big](1-q^n) + n pq^n \big)\nonumber \\ 
&\overset{(a)}{<} 0, \text{ if }  r-1- p(n+(r/p^*)-1)<-1, 
\label{eq:u_r_der_B_2} 
\end{align}
where $(a)$ follows from the fact that $npq^n - (1-q^n) < 0$ for $p \in (0,1)$.  So, $u_r(p,n,g)$ decreases for $r/(n+(r/p^*)-1)<p<p^*$ implying that
$$p^*_r\leq \frac{r}{n+(r/p^*)-1}.$$
Since $g(p)$ is an increasing function for $p \in (0,p^*)$, we get
\begin{equation}
c(r)=\frac{g(p^*_r)}{r-1} \leq \frac{1}{r-1} g\left(\frac{r}{n+(r/p^*)-1} \right), r\geq 3. 
\label{eq:c_ub_r_B}
\end{equation}
In the next lemma, we present upper bounds on $\frac{1}{r-1} \ g\left(\frac{r }{n+r/p^*-1}\right)$ that are independent of $r$ and are in the form of \eqref{eq:rate_of_fall}.
\begin{lemma} \label{lem:type-B-c}
    \begin{enumerate}
    \item  If either $g(p)/p$ is non-increasing, or $n < 1+\frac{4}{(1-p^*)^2}$,    
    $$c(r) \leq 0.5 \ g\left( \frac{3}{n+({3}/{p^*})-1}\right).$$ To get the form of \eqref{eq:rate_of_fall}, set $b_0 = b_3 = 0.5, b_1 = b_4 = 3 $ and $ b_2 = b_5 = ({3}/{p^*})-1$. \label{lem-part:g(p)/p_decreasing}

    \item  If $n \geq 1+\frac{4}{(1-p^*)^2}$,  
    $$c(r) \leq \max \left( 0.5 \ g\left( \frac{3}{n+(3/p^*)-1}\right), \frac{1}{r_4-1}\ g\left(\frac{r_4}{n+(r_4/p^*)-1}\right)\right),$$ 
    where $r_4 = 0.5\ (n-1) \ (1-p^*)\ \left(1 + \sqrt{1 - \frac{4}{(n-1) (1-p^*)^2}} \right)$. The above bound is directly in the form of \eqref{eq:rate_of_fall}.
    \label{lem-part:disc>0}
    \end{enumerate}
\end{lemma}
\begin{proof}
A proof is given in Appendix \ref{subsec:type-B-c_proof}.
\end{proof}

\section{Proof of Corollaries}
\label{sec:corollary_proof}
We begin with deriving upper bounds for $u^*_2(n,g(p))$ for  $g(p)=p^{\alpha}$ and $g(p)=p\log(1/p)$. These are used in the simplifications of the tail bound expression.
\subsection{Upper bounds on $u^*_{2}(n,g(p))$} 
\label{subsec: Var_bound} 
Recall that for $g(p) = p$, $u_{2}(n,p) = p (1-p)^n (1-(1-p)^n)$, $u^*_{2}(n,p) = \max_{p \in (0,1) }\,u_2(n,p)$ and
\begin{align*}
 \gamma = \max_{np}\, np\, e^{-np} (1-e^{-np}),\,p\in(0,1).
\end{align*}
Since $\argmax p(1-p)^n = 1/(n+1)$ and $1-(1-p)^n$ is increasing in $p$, we have that $\argmax_{0<p<1} u_{2}(n,p) > \frac{1}{n+1}$.

If $p \geq \ln 2/n$ or $e^{-np} \leq 1/2$, since $0\le 1-p\le e^{-p}$, we have $$e^{-np}+(1-p)^n \leq 1 \text{ if }p \geq \ln 2/n.$$ 
Multiplying by the positive quantity $np(e^{-np}-(1-p)^n)$ and rearranging, we get
\begin{equation}
    np (1-p)^n (1-(1-p)^n) \leq np e^{-np} (1-e^{-np}) \text{ if } p\geq\frac{\ln 2}{n}. \label{eq:u_2^*_bound}
\end{equation}
Since $\argmax_{0<p<1} u_{2}(n,p) >\frac{1}{n+1} > \frac{\ln 2}{n}$ for $n \geq 3$, we have $u^*_{2}(n,p) \leq \gamma /n$, for $n \geq 3$. 

The generalization for $\alpha\ge 1$ is very similar to the proof above. So, we omit the details and present only the final result. For $\alpha\ge 1$, 
$$u^*_{2}(n,p^{\alpha}) = \max_{p \in (0,1)} \,p^{2\alpha-1} (1-p)^n (1-(1-p)^n).$$
For $n \geq (2\alpha-1)\ln 2/(2\alpha-1 - \ln 2)$, the following bound holds:
$$u_2^*(n,p^{\alpha}) \leq \gamma_{\alpha}/n^{2\alpha-1},$$ 
where $\gamma_{\alpha} = \max_{t > 0} t^{2\alpha-1} e^{-t}(1-e^{-t})$. 

For the case of missing Shannon entropy $H_0(p),$ we have 
\begin{align}
    u^*_{2}(n,p\log_2(1/p)) \ & = \  \max_{p \in (0,1) }\,u_2(n,p\log_2(1/p)) \nonumber \\ 
    & = \  \max_{p \in (0,1) }\, p \  (\log_2 (1/p) )^2 \ (1-p)^n (1-(1-p)^n) \nonumber \\ 
   & \overset{(a)}{=} \ \max_{p \in (0,1) }\, p \ (\log_2 (np) - \log_2 n )^2 \ (1-p)^n (1-(1-p)^n) \nonumber \\ 
    & \overset{(b)}{\leq} \ ( \log_2 n )^2 \ u^*_{2}(n,p)  \ \leq \ ( \log_2 n )^2 \ \gamma/n, \ n \geq 3
\end{align}
where we get $(a)$ by using $\log_2 p = \log_2 (np/n) = \log_2(np) - \log_2 n$ and $(b)$ by using $\log_2 n \geq \log_2 (np).$ Therefore, $ u^*_{2}(n,p\log_2(1/p)) \leq ( \log_2 n )^2 \ \gamma/n,$  $n \geq 3.$ 

The upper bound on $u^*_2(n,p)$ can be extended to an upper bound on $u^*_2(n,g(p))$ in cases where $g(p)/p$ is bounded by using the following:
\begin{align}
    u^*_2(n,g(p)) &= \max_{0\le p\le 1} (g(p)^2/p)(1-p)^n (1-(1-p)^n)\nonumber\\
    &\le \max_{0\le p\le 1} (g(p)/p)^2 \max_{0\le p\le 1} p(1-p)^n(1-(1-p)^n)\nonumber\\
    &\le (\gamma/n) \max_{0\le p\le 1} (g(p)/p)^2.\quad (n>3)
\end{align}


\subsection{Proof of Corollary \ref{cor:left_tail_special_case}}
 Using Part 1 of Theorem \ref{thm:con_G0} with $u_2^*(n,p^{\alpha}) \leq \gamma_{\alpha}/n^{2\alpha-1}$ and $u_2^*(n,p \log_2(1/p)) \leq (\log_2 n)^2 \gamma/n$ (as shown above) in the Chernoff method results in the left tail bounds presented in Corollary \ref{cor:left_tail_special_case} for  $M_{0,\alpha}, \alpha \geq 1$ and $H_{0}$, respectively.
    
\subsection{Proof of Corollary \ref{cor:right_tail_special_case}}
The proof for each part of Corollary \ref{cor:right_tail_special_case} is given below. For every $g(p)$, we identify the scale parameter $c$ either from Lemma \ref{lem:type-A-c} if $g(p)$ is of Type A, or from Lemma \ref{lem:type-B-c} if $g(p)$ is of Type B. We bound the scale parameter, if necessary. Then, the variance parameter is found by bounding $u^*_2(n,g)$ using one of the methods described in Section \ref{subsec: Var_bound}. Finally, $c$ and $v$ are used in the strongly sub-Gamma right tail bound \eqref{eq:strongly-sub-gamma-tail-bound}.
\begin{enumerate}
    \item For $g(p)=p$, we get $c = 3/(2(n+2))$ from Lemma \ref{lem:type-A-c} Part \ref{lem-part:mu < 1}. \textcolor{blue}{Using $c = 3/(2(n+2))$ and $v = u^*_{2}(n,p)\le \gamma/n$ (as shown above) in \eqref{eq:strongly-sub-gamma-tail-bound}, \eqref{eq:strongly-sub-gamma-tail_ub} gives \eqref{eq:right_tail_derived_ssg_M0}, \eqref{eq:right_tail_derived_ssg_M0_ub} respectively.} 
    \textcolor{blue}{Using $v = (1+2/n) \ \psi/3n, c = 3/2(n+2),$ and $ a = [a_2] =  [\gamma/n - v]$ gives \eqref{eq:right_tail_derived_poly_ssg_M0_ub}.}
    
    \item For $g(p)=p^{\alpha}$, $\alpha > 1$ and $n$ sufficiently large, using Lemma \ref{lem:type-A-c}, Part \ref{lem-part:mu>1_disc>0} with $\mu=\alpha$, it can be shown that 
    \begin{align}
        c &= \frac{1}{{r_2} -1 }({\alpha{r_2}}/{(n-1+\alpha {r_2})})^{\alpha}\le \frac{2\alpha(2(\alpha-1)/(\alpha+1))^{\alpha}}              {(n-1)(\alpha-1)-2\alpha},
    \end{align}
where $r_2$ is as defined in Lemma \ref{lem:type-A-c}, and the inequality is obtained by using $0\le \sqrt{1-4\alpha^2/((n-1)(\alpha-1)^2)}\le 1$ appropriately. Using the $c$ above in \eqref{eq:strongly-sub-gamma-tail-bound} along with $v = u_2^*(n,p^{\alpha}) \leq \gamma_{\alpha}/n^{2\alpha -1}$ (as shown above) gives  \eqref{eq:right_tail_derived_ssg_M0_alpha}.
    
\item For $g(p) = p \log_2(1/p)$, from Lemma \ref{lem:type-B-c} Part \ref{lem-part:g(p)/p_decreasing}, we get $c = 0.5\ {(\log_2 n_0)} / {n_0 }$, where $n_0 = (1/3)(n-1) + e$. Using this value of $c$ in \eqref{eq:strongly-sub-gamma-tail-bound} along with $v=u_2^*(n,p \log_2(1/p) ) \leq  (\log_2 n)^2 \gamma /n$ (as shown above) gives \eqref{eq:right_tail_derived_ssg_H_0}. Using the taylor series bound $\ln(1+x) \leq x - 0.5x^2 + (1/3)x^3, x \geq 0,$ with $x = \epsilon \log_2(n_0)/(\gamma_n)$ gives  \eqref{eq:right_tail_derived_ssg_H_0_ub}.
\end{enumerate}

\section{Conclusion and Future Directions}
\label{sec:conclusion}
We have generalized the notion of missing mass to missing $g-$ mass of a function $g : [0,1] : \to \mathbb{R}$ \textcolor{blue}{motivated by the use of special cases like missing mass of order$-\alpha$ and missing Shannon entropy in estimating the closeness of the missing probabilities to uniformity.} Estimation and concentration of missing mass of functions was studied and several new results were shown. \textcolor{blue}{The estimates of the missing mass of functions were applied to estimate the missing probabilities together with their multiplicities and test the closeness of the missing probabilities to uniformity.}
In particular, by generalizing Good-Turing estimators, we showed that estimation better than that of a plugin estimator is possible for several interesting classes of functions. However, there are important classes of functions where linear estimators perform only as good as the plugin estimator. Non-linear estimators are an interesting possibility to explore in future work.

As far as concentration results are concerned, we introduced two new notions of concentration, named strongly sub-Gamma and filtered sub-Gaussian, that result in tail bounds better than that of sub-Gaussian concentration for missing mass. In other situations where sub-Gaussian concentration has a sub-optimal variance factor, these new notions, particularly the idea of filtering a Gaussian distribution, is worth exploring.


\bibliographystyle{IEEEtran}
\bibliography{main}

\newpage
\appendix
\subsection{Proof of Lemma \ref{lem:DirichletPrior}}
\label{appsec:dirichletproof}
Recall that $P = (p_1, p_2, \ldots, p_k)$ on $\Delta_k$ has a Dirichlet distribution with parameters $\bbeta = (\beta_1, \beta_2, \ldots, \beta_k)$.
Let $\beta_0 \triangleq \sum_{i=1}^k \beta_i$. Let $\Gamma(u)\triangleq\int_0^{\infty} x^{u-1}e^{-x}dx$ denote the Gamma function, and let $\tau(u,v) \triangleq {\Gamma(u+v)}/{\Gamma(u)}$. Recall the notation $N_x(X^n)=\sum_{i=1}^n I(X_i=x)$ denoting the number of occurrences of $x$ in $X^n$. The following properties of the Dirichlet distribution are useful:
\begin{align}
    E[p_i^a] &= \tau(\beta_i,a)/\tau(\beta_0,a),\ a>0\label{eq:expectation1}\\
    E[p_i^a p_j^b] &= \tau(\beta_i,a) \tau(\beta_j, b)/\tau(\beta_0,a+b),\ a,b>0\label{eq:expectation2}\\
    P\vert X^n &\sim \text{Dirichlet} (\beta_1+F_1(X^n),\beta_2+F_2(X^n),\ldots,\beta_k+F_k(X^n)).\label{eq:conjugate}
\end{align}
Since $M_{0,\alpha} = \sum_{x \in \cX} p_x^{\alpha} I(N_x(X^n) =0)$ and $N_x$ is a function of only $X^n$, we get
$$E\left[M_{0,\alpha}(X^{n},P)\vert X^{n}\right] = \sum_{x \in \mathcal{X}} E[p_x^{\alpha}|X^n] \ I(N_x = 0).$$ 
The conditional variance of $M_{0,\alpha}(X^{n},\ P)$ simplifies as follows:
\begin{align}
&\text{Var}\left[M_{0,\alpha}(X^{n},\ P)\vert X^{n}\right]  = E\left[\left( \sum_{x \in \mathcal{X}} I(N_x = 0) (p_x^{\alpha} - E[p_x^{\alpha}|X^n]) \right)^2\bigg\rvert X^n \right] \nonumber \\ 
&= \sum_{x,y \in \mathcal{X} : x \neq y} I(N_x = N_y = 0) \big(E[p_x^\alpha p_y^\alpha|X^n]-E[p_x^\alpha|X^n]E[p_y^\alpha|X^n]\big)\nonumber\\
&\qquad\qquad + \sum_{x \in \mathcal{X}} I(N_x = 0) \big(E[p_x^{2\alpha}|X^n]-E[p_x^{\alpha}|X^n]^2\big)\nonumber\\
&\overset{(a)}{=}  \sum_{x,y \in \mathcal{X} : x \neq y} I(N_x = N_y = 0)  \ { \tau( \beta_x +N_x,\alpha)}    \tau( \beta_y +N_y,\alpha) \left( \frac{1 }{\tau(\beta_0 +n,2\alpha) }-  \frac{1 }{\tau^2( \beta_0 +n,\alpha) } \right)  \nonumber \\ 
& \qquad + \sum_{x \in \mathcal{X}} I(N_x = 0) \left( \frac{ \tau(\beta_x+N_x, 2\alpha ) }{ \tau(\beta_0+n, 2\alpha ) }  -  \frac{ \tau^2(\beta_x+N_x, \alpha ) }{ \tau^2(\beta_0+n, \alpha ) }  \right) \nonumber\\
&\overset{(b)}{=}  \sum_{x,y \in \mathcal{X} : x \neq y} I(N_x = N_y = 0)  \ { \tau( \beta_x,\alpha)}    \tau( \beta_y,\alpha) \left( \frac{1 }{\tau(\beta_0 +n,2\alpha) }-  \frac{1 }{\tau^2( \beta_0 +n,\alpha) } \right)  \nonumber \\ 
& \qquad + \sum_{x \in \mathcal{X}} I(N_x = 0) \left( \frac{ \tau(\beta_x, 2\alpha ) }{ \tau(\beta_0+n, 2\alpha ) }  -  \frac{ \tau^2(\beta_x, \alpha ) }{ \tau^2(\beta_0+n, \alpha ) }  \right), \label{eq:conditionalvariance}
\end{align}
where $(a)$ follows by using \eqref{eq:conjugate} and \eqref{eq:expectation1}, \eqref{eq:expectation2}, and $(b)$ follows because of the presence of the indicators $I(N_x = N_y = 0)$ and $I(N_x = 0)$. Taking expectation over $X^n$ on both sides of \eqref{eq:conditionalvariance},
\begin{align}
   &E_{X^n}\left[\text{Var}\left[M_{0,\alpha}(X^{n},P)\vert X^{n}\right]\right]\nonumber\\ 
   &= \sum_{x,y \in \mathcal{X} : x \neq y} \pr(N_x = N_y = 0)  \ { \tau( \beta_x,\alpha)}  \tau( \beta_y ,\alpha) \left( \frac{1 }{\tau(\beta_0 +n,2\alpha) }-  \frac{1 }{\tau^2( \beta_0 +n,\alpha) } \right) \nonumber \\ 
   & \qquad + \sum_{x \in \mathcal{X}} \pr(N_x = 0) \left( \frac{ \tau(\beta_x, 2\alpha ) }{ \tau(\beta_0+n, 2\alpha ) }  -  \frac{ \tau^2(\beta_x, \alpha ) }{ \tau^2(\beta_0+n, \alpha ) }  \right)\nonumber\\
   &\overset{(a)}{=} \sum_{x \in \mathcal{X}} \frac{\tau(\beta_0-\beta_x, n)}{\tau(\beta_0,n)} \left( \frac{ \tau(\beta_x, 2\alpha ) }{ \tau(\beta_0+n, 2\alpha ) }  -  \frac{ \tau^2(\beta_x, \alpha ) }{ \tau^2(\beta_0+n, \alpha ) }  \right)\nonumber \\ 
   & \qquad + \sum_{x,y \in \mathcal{X} : x \neq y} \frac{\tau(\beta_0-\beta_x-\beta_y, n)}{\tau(\beta_0,n)}\, \tau( \beta_x,\alpha)\tau( \beta_y ,\alpha) \left( \frac{1 }{\tau(\beta_0 +n,2\alpha) }-  \frac{1 }{\tau^2( \beta_0 +n,\alpha) } \right) ,
   \label{eq:lb_beta}
\end{align}
where $(a)$ follows because  
\begin{align*}
    \pr(N_x(X^n)=0)&=E[(1-p_x)^n]=\tau(\beta_0-\beta_x,n)/\tau(\beta_0,n),\\
    \pr(N_x(X^n)=N_y(X^n)=0)&=E[(1-p_x-p_y)^n]= \tau(\beta_0-\beta_x-\beta_y, n)/\tau(\beta_0,n).
\end{align*}
Let $D_{n,c}$ denote the Dirichlet distribution with $\beta_i = 1/n$, $i = 1,2,\ldots,k$ and $k = cn^2$ with $c > 0$. Let 
$$T(D_{n,c})\triangleq E_{X^n}\left[\text{Var}\left[M_{0,\alpha}(X^{n},D_{n,c})\vert X^n\right]\right].$$
Setting $P=D_{n,c}$ in \eqref{eq:lb_beta}, we get
\begin{align}
T(D_{n,c})&= cn^2\,\frac{\tau(cn-1/n, n)}{\tau(cn,n)} 
\bigg(\frac{\tau(1/n,2\alpha)}{\tau((c+1)n,2\alpha)}-\frac{\tau^2(1/n,\alpha)}{\tau^2((c+1)n,\alpha)}\bigg)\nonumber\\ 
&\ + cn^2(cn^2-1)\,\frac{\tau(cn-2/n,n)}{\tau(cn,n)}\,\tau^2(1/n,\alpha)
\left(\frac{1}{\tau((c+1)n,2\alpha)}-\frac{1}{\tau^2((c+1)n,\alpha)}\right),\nonumber\\
&\triangleq cn^2A_1(A_2 - A_3) + cn^2(cn^2-1)A_4 A_5(A_6 - A_7),
\label{eq:lb_beta_(1/n)} 
\end{align} 
where $A_i$, $1\le i\le 7$, denote the corresponding terms in the previous equation. We consider two different cases.

We use
$$\tau(u,v) = \Gamma(u+v)/\Gamma(u) = \prod_{l= 0}^{v-1} (u+l),\ v \in \mathbb{N},$$ 
to simplify $A_i$ and find their dominant terms as follows.
\begin{align}
    &A_1 = \prod_{l = 0}^{n-1} \left( 1 - \frac{1}{n(cn+l)}\right) =_n 1,\ A_2 = \prod_{l=0}^{2\alpha-1}\frac{1/n + l}{(c+1)n+l} =_n \frac{(2\alpha-1)!}{(c+1)^{2\alpha}n^{2\alpha+1}},\nonumber\\
    &A_3 = \prod_{l=0}^{\alpha-1}\frac{(1/n + l)^2}{((c+1)n+l)^2} =_n O(1/n^{2\alpha+2}),\ A_4 = \prod_{l=0}^{n-1}\left(1-\frac{2}{n(cn+l)}\right)=_n 1,\nonumber\\ 
    &A_5 = \prod_{l=0}^{\alpha-1} {(1/n + l)^2 } =_n \frac{((\alpha-1)!)^2}{n^2},\nonumber\\
    &A_6 = \prod_{l=0}^{2\alpha-1}\frac{1 }{((c+1)n+l) } =_n \frac{1}{((c+1)n)^{2\alpha}+(\alpha-1)(2\alpha-1)((c+1)n)^{2\alpha-1}},\nonumber\\ 
    &A_7 = \prod_{l=0}^{\alpha-1} \frac{1 }{( (c+1)n+l)^2 } =_n \frac{1}{((c+1)n)^{2\alpha}+(\alpha-1)(\alpha-2)((c+1)n)^{2\alpha-1}}.\nonumber
\end{align}
For $A_6$ and $A_7$, the first two terms are retained because they are being subtracted and the leading terms are identical. Using the above in \eqref{eq:lb_beta_(1/n)}  we get 
   \begin{align}
    &T(D_{n,c}) =_n  cn^2 \left(\frac{(2\alpha-1)!}{(c+1)^{2\alpha}n^{2\alpha+1}} - O(1/n^{2\alpha+2})\right)+ \ cn^2(cn^2-1) \frac{((\alpha-1)!)^2}{n^2}  \nonumber \\ 
    &\qquad\qquad\qquad \quad \bigg( \frac{1}{((c+1)n)^{2\alpha}+(\alpha-1)(2\alpha-1)((c+1)n)^{2\alpha-1}} \nonumber \\ & \qquad \qquad \qquad \qquad \qquad  -\frac{1}{((c+1)n)^{2\alpha}+(\alpha-1)(\alpha-2)((c+1)n)^{2\alpha-1}}\bigg)  \nonumber \\
    &\qquad\qquad \geq_n \frac{c(2\alpha-1)!}{(c+1)^{2\alpha}n^{2\alpha-1}}- \ cn^2(cn^2-1) \frac{((\alpha-1)!)^2}{n^2}\frac{\alpha^2}{((c+1)n)^{2\alpha+1}}\nonumber\\
    &\qquad\qquad=_n \frac{c/(c+1)}{((c+1)n)^{2\alpha -1}} \left( (2\alpha -1)! - \frac{c(\alpha!)^2}{c+1}\right).\label{eq:lb_int}
\end{align} 
This completes the proof of Lemma \ref{lem:DirichletPrior}.

\subsection{Proof of Lemma \ref{lem:type-A-c} } \label{subsec:type-A-c_proof}
Let $p(r)\triangleq r\mu/(n+r\mu-1)$ and $h(p)\triangleq g(p)/p$. The general idea of the proof is to start with \eqref{eq:c_ub_r_A}, which states 
$$c(r)\le c_{\max}(r)\triangleq g(p(r))/(r-1),$$ 
and try to maximize $c_{\max}(r)$ over $r\ge 3$.

\subsubsection{Case 1(a) ($0<\mu<1$)}

Since $g'(p)=p h'(p)+h(p)\le \mu h(p)$ or $h'(p)\le (\mu-1)h(p)/p < 0$, we see that $h(p)$ decreases with $p$. Rewriting $c_{\max}(r)$ in terms of $h(p(r))$, we have
$$c(r) \leq \frac{p(r)}{r-1} h(p(r)) \overset{(a)}{\leq} \frac{p(3)}{2}\,h(p(r)) \overset{(b)}{\le} \frac{p(3)}{2}\,h(p(3)) = 0.5g(p(3)),$$ 
where $(a)$ follows because $p(r)/(r-1)$ decreases with $r$ and $r\ge 3$, and $(b)$ follows because $h(p(r))$ decreases with $r$ ($p(r)$ increases with $r$, $h(p)$ decreases with $p$) and $r\ge 3$.

This completes the proof for Case 1(a).

\subsubsection{Case 1(b) ($\mu>1$ and $n < 1 + 4(\mu/(\mu-1))^2$)}

The derivative of $c_{\max}(r)$ can be simplified as follows.
\begin{align}
c'_{\max}(r)&=\frac{1}{r-1} \bigg(\frac{(n-1)\mu}{(n+r\mu-1)^2}\,g'(p(r)) - \frac{1}{(r-1)}\,g(p(r))\bigg)\nonumber\\
&\overset{(a)}{\le} g(p(r))\,\frac{1}{r-1}\bigg(\frac{(n-1)\mu}{r(n+r\mu-1)}
              - \frac{1}{r-1}\bigg)\nonumber\\
&= -\mu\, g(p(r))\,\frac{r^2-(n-1)(1-1/\mu)r+(n-1)}{r(r-1)^2(n+r\mu-1)}\nonumber\\
&= -\mu\, g(p(r))\,\frac{\big[r-(n-1)(1-1/\mu)/2\big]^2+(n-1)\big[1-(n-1)(1-1/\mu)^2/4\big]}{r(r-1)^2(n+r\mu-1)},\label{eq:cmax_derivative}
\end{align}
where $(a)$ uses $g'(p(r))\le \mu\,g(p(r))/p(r)$. 

If $1-(n-1)(1-1/\mu)^2/4>0$ or $n < 1 + 4\mu^2/(\mu-1)^2$, we see that $c'_{\max}<0$ and $c_{\max}(r)$ decreases with $r$. So, for $r\ge 3$,
$$c(r)\le c_{\max}(r) \le c_{\max}(3) = 0.5g(p(3)).$$
This completes the proof for Case 1(b).

\subsubsection{Case 2 ($\mu>1$ and $n \ge 1 + 4(\mu/(\mu-1))^2$)}

From \eqref{eq:cmax_derivative}, $c'_{\max}(r)$ is negative for $r>r_2$, where $$r_2\triangleq 0.5(n-1)(1-1/\mu)(1 + \sqrt{1-4\mu^2/((n-1)(\mu -1)^2)})$$ 
is the largest of the roots of the quadratic polynomial in the numerator of \eqref{eq:cmax_derivative}. So, we have
$$c(r)\le c_{\max}(r) \le \max \bigg(c_{\max}(3), c_{\max}(r_2)\bigg) = \max \bigg(g(p(3))/2, g(p(r_2))/(r_2-1) \bigg).$$
This completes the proof for Case 2.

\subsection{Proof of Lemma \ref{lem:type-B-c} } \label{subsec:type-B-c_proof}
The proof mirrors the above proof of Lemma \ref{lem:type-A-c} with some minor changes to the functional forms, and we skip the details.

\end{document}